\documentclass[a4paper,10pt]{amsart}
\usepackage[english]{babel}
\usepackage{amsmath,tikz-cd}
\usepackage{amssymb}
\usepackage[T1]{fontenc}
\usepackage[utf8]{inputenc}
\usepackage{lmodern}
\usepackage{mathrsfs}
\usepackage{enumerate}
\usepackage[shortlabels]{enumitem}
\usepackage{mathtools,yfonts}
\usepackage{tikz,tkz-euclide}
\usepackage{booktabs}
\usetikzlibrary{arrows,arrows.meta,decorations.markings}
\usepackage{makecell}
\usepackage{xcolor}
\definecolor{DarkRed}{RGB}{173,0,0}
\definecolor{LightRed}{RGB}{201,0,0}
\usepackage[
    colorlinks=true,
    linkcolor=DarkRed,
    urlcolor=LightRed,
    citecolor=LightRed
]{hyperref}
\usepackage{tcolorbox}
\usepackage[]{algorithm2e}

\newtheorem{thm}{Theorem}[section]
\newtheorem{Lemma}[thm]{Lemma}
\newtheorem{Proposition}[thm]{Proposition}
\newtheorem{Corollary}[thm]{Corollary}

\newtheorem{teo}[thm]{Theorem}

\theoremstyle{definition}

\newtheorem{Definition}[thm]{Definition}
\newtheorem{Remark}[thm]{Remark}

\newtheorem{oss}[thm]{Remark}

\definecolor{wwwwww}{rgb}{0.4,0.4,0.4}

\hypersetup{pdfpagemode=UseNone}
\hypersetup{pdfstartview=FitH}

\begin{document}

	\title{On the arithmetic of rational hypersurfaces in toric varieties}

	\author[Gianluca Grassi]{Gianluca Grassi}
	\address{\sc Gianluca Grassi\\ Dipartimento di Matematica e Informatica, Università di Ferrara, Via Machiavelli 30, 44121 Ferrara, Italy}
	\email{gianluca.grassi@unife.it}

	\date{\today}
	\subjclass[2020]{Primary 12F20, 12E10, 14E08, 14M20, 12F10; Secondary 14G05.}
	\keywords{Transcendental field extensions, rationality, rational points, hypersurfaces.}
	
	\maketitle
\begin{abstract}
We study a family of rational hypersurfaces of odd degree in $\mathbb P^{2n+1}$ with prescribed multiplicity along two skew conjugate $n$-planes. After blowing up these two linear spaces, their strict transforms are realized as divisors $\widetilde X^{2n}\subset \mathcal T$ of multidegree $(2d+1,-d,-d)$ in a smooth toric variety $\mathcal T$. Over the quadratic extension $k(\xi)/k$, with $\xi^2=-3$, the Cox ring of $\mathcal T$ is $k(\xi)[z_0,\ldots,z_{2n+1},w_+,w_-]$, graded by $\deg(z_{2i})=(1,-1,0)$, $\deg(z_{2i+1})=(1,0,-1)$, $\deg(w_+)=(0,1,0)$ and $\deg(w_-)=(0,0,1)$. We prove that $\widetilde X^{2n}$ is $k$-rational, and derive results on the distribution of its rational points over number fields and finite fields. The case $d=1$ recovers the even-dimensional Fermat cubic.
\end{abstract}

	\tableofcontents
	
\section{Introduction}

Throughout the paper, unless explicitly stated otherwise, $k$ denotes a field, $\bar k$ a separable closure, and $d\geq 1$, $n\geq 1$ are integers. In the sections based on the quadratic extension $k(\xi)/k$ we assume $\operatorname{char}(k)\neq 2,3$ and $\xi^2=-3$. The two $n$-planes used in the construction are defined over $k(\xi)$ and are exchanged by the non-trivial element of $\operatorname{Gal}(k(\xi)/k)$. In the finite-field statements the condition $q\equiv 5\pmod 6$ implies that $-3$ is not a square in $\mathbb F_q$, hence these planes are not individually defined over the base field. Characteristic $2$ is treated separately in Section~\ref{sec:char-two}. We do not use the restriction-of-scalars construction in characteristic $3$.

Fix an integer $n\geq 1$, a field $k$ with quadratic extension $k(\xi)$, where $\xi^2=-3$, and an integer $d\geq 1$. Let $R=k(\xi)[z_0,\ldots,z_{2n+1},w_+,w_-]$ be the Cox ring with $\mathbb Z^3$-grading given by $\deg(z_{2i})=(1,-1,0)$, $\deg(z_{2i+1})=(1,0,-1)$, $\deg(w_+)=(0,1,0)$ and $\deg(w_-)=(0,0,1)$, for $i=0,\ldots,n$. Let $\mathfrak B$ be the irrelevant ideal
$$
\mathfrak B
=
(z_0,\ldots,z_{2n+1})
\cap
(w_+,z_1,z_3,\ldots,z_{2n+1})
\cap
(w_-,z_0,z_2,\ldots,z_{2n}).
$$
We denote by $\mathcal T=(\operatorname{Spec}R\setminus V(\mathfrak B))/(\mathbb G_m)^3$ the resulting smooth toric variety. Over $k(\xi)$ this is the blow-up of $\mathbb P^{2n+1}$ along the two coordinate $n$-planes $\{z_0=z_2=\cdots=z_{2n}=0\}$ and $\{z_1=z_3=\cdots=z_{2n+1}=0\}$; over $k$ it descends by Galois symmetry.

Define
$$
\widehat S
=
\sum_{i=0}^n z_{2i}^d z_{2i+1}^d(w_+z_{2i}+w_-z_{2i+1}),
\qquad
\widehat D
=
\sum_{i=0}^n z_{2i}^d z_{2i+1}^d(w_+z_{2i}-w_-z_{2i+1}).
$$
Both polynomials have multidegree $(2d+1,-d,-d)$. Set $\widehat F:=3\widehat D-\xi\widehat S$, and let $\widetilde X^{2n}:=V(\widehat F)\subset \mathcal T_{k(\xi)}$. The equation $\widehat F=0$ is written over $k(\xi)$; the corresponding $k$-model is obtained by the Galois-symmetric construction described below. In particular, on $k$-rational points the condition $\widehat F=0$ may be read through its conjugation-invariant real and imaginary parts.

In $\mathbb P^{2n+1}$ consider the odd-degree hypersurface
\begin{equation}
\label{X2n_equation}
X^{2n}
=
\left\{
\sum_{i=0}^n
(x_{2i}+x_{2i+1})
(x_{2i}^2-x_{2i}x_{2i+1}+x_{2i+1}^2)^d
=
0
\right\}
\subset \mathbb P^{2n+1}.
\end{equation}
It has degree $2d+1$ and multiplicity $d$ along two skew conjugate $n$-planes. Taking lines through conjugate points of multiplicity $d$ gives a rational map $\overline{\varphi}_{2n}\colon \mathbb P^{2n}\dashrightarrow X^{2n}$; see the formulas in \eqref{projectivepolynomialsAB}--\eqref{projectivephis}.

\begin{teo}
\label{thm:main-toric}
With the notation above, the following hold.
\begin{enumerate}[label=\textup{(\roman*)}]
\item The strict transform $\widetilde X^{2n}\subset \mathcal T$ cut out by $\widehat F=0$ has class $(2d+1,-d,-d)$ and is birational to $\mathbb P^{2n}$. That is, the rational map $\overline{\varphi}_{2n}\colon\mathbb P^{2n}\dashrightarrow X^{2n}$ factors through $\mathcal T$ and induces a birational correspondence $\widetilde X^{2n}\dashrightarrow X^{2n}$. For $d=1$ this recovers the even-dimensional Fermat cubic of \cite{msslxa}.

\item If $k=\mathbb F_q$, $q\equiv 5\pmod 6$, and $\gcd(2d+1,q-1)=1$, then
$$
\big|\widetilde X^{2n}(\mathbb F_q)\big|
=
\big|\mathbb P^{2n}(\mathbb F_q)\big|
=
\frac{q^{2n+1}-1}{q-1}.
$$
For $n=1$ an explicit formula holds for all $q$; see Corollary~\ref{numberpoints}.

\item For any number field $K$ and the height $H$ induced by $\mathcal O_{X^{2n}}(1)$, one has
$$
\#\{x\in \widetilde X^{2n}(K): H(x)\leq B\}
\gg
B^{\frac{2n+1}{2d+2}}
\qquad
(B\to\infty),
$$
with the improved lower bound $B^{\frac{3}{2d+1}}$ when $n=1$.
\end{enumerate}
\end{teo}

\begin{oss}
Over $k(\xi)$ the ambient variety $\mathcal T$ is a smooth split toric variety. The Manin--Peyre prediction and asymptotics for toric varieties are established in \cite{BT98}; see also \cite{MS20} on approximation, \cite{Pie24} for the hyperbola method, and \cite{PS24,Bon24} for bounded-height points on toric subvarieties and multi-height refinements. Our contribution is a concrete toric realization of the strict transform $\widetilde X^{2n}$ via a single Cox equation of multidegree $(2d+1,-d,-d)$ over $k(\xi)$, together with the explicit birational bridge to $X^{2n}$ yielding finite-field identities; see Corollary~\ref{numberpoints}, Proposition~\ref{pointcountX2n} and Proposition~\ref{pointcountY2n-2}.
\end{oss}

\subsubsection*{Structure of the paper.}
We begin with the Weil restriction reformulation and the explicit parametrization $\overline{\varphi}_{2n}$ in affine and projective form; see \eqref{projectivepolynomialsAB}--\eqref{projectivephis}. The characteristic $2$ case is handled by a direct Cox-style parametrization. We then introduce a quadro-cubic Cremona transformation relating the two parametrizations and derive point-counting identities over finite fields; see Corollary~\ref{numberpoints} and Proposition~\ref{pointcountX2n}. Finally, we prove height bounds and describe the Cox model of $\mathcal T$. In this framework $\widetilde X^{2n}$ is given by the single equation $\widehat F=0$ over $k(\xi)$, while its descent to $k$ is encoded by the Galois-symmetric construction.

\section{Rationality construction via restriction of scalars}
\label{sec:restriction-scalars}

We begin with a general rationality construction for hypersurfaces of odd degree in $\mathbb P^{2n+1}$.

\begin{Proposition}
\label{grassmannian}
Let $X\subset \mathbb P^{2n+1}$ be a hypersurface of degree $2d+1$ over a field $k$, with $\operatorname{char}(k)\neq 2$. Assume that, over a quadratic extension of $k$, the hypersurface $X$ contains two skew conjugate $n$-planes $H_+$ and $H_-$ with multiplicity $d$. Then $X$ is rational over $k$.
\end{Proposition}

\begin{proof}
Let $\Lambda\subset \mathbb P^{2n+1}$ be a $(n-1)$-plane defined over $k$ such that $\Lambda\cap(H_+\cup H_-)=\emptyset$. A general $(n+1)$-plane $H\subset \mathbb P^{2n+1}$ defined over $k$ and containing $\Lambda$ intersects $H_+$ and $H_-$ in two conjugate points $p_+$ and $p_-$. Since $p_+$ and $p_-$ have multiplicity $d$ on $X$, the line $L_H=\langle p_+,p_-\rangle$ meets $X$ in a residual point $x_H$. Since $L_H$ and $X$ are defined over $k$ and the two length-$d$ intersections at $p_+$ and $p_-$ are conjugate, the residual point $x_H$ is defined over $k$.

The $(n+1)$-planes containing $\Lambda$ are parametrized by $\mathbb G(1,n+1)$, which has dimension $2n$ and is rational over $k$. We therefore obtain a rational map $\phi\colon \mathbb G(1,n+1)\dashrightarrow X$, $H\mapsto x_H$. This map is generically injective. Indeed, for a general point $x\in X$, there is a unique line through $x$ meeting the two skew complementary $n$-planes $H_+$ and $H_-$; it is obtained from the decomposition of $x$ with respect to their linear span. Hence the points $p_+$ and $p_-$, and therefore the plane $\langle \Lambda,p_+,p_-\rangle$, are recovered from $x$. Thus $\phi$ is birational, and $X$ is rational over $k$.
\end{proof}

We now reformulate the construction in affine space via Weil restriction. Let $L/k$ be a finite field extension, and let $X=\operatorname{Spec}(L[x_1,\ldots,x_n]/(f_1,\ldots,f_m))$ be an affine variety over $L$. Fix a $k$-basis $e_1,\ldots,e_s$ of $L$, introduce variables $y_{i,j}$ for $i=1,\ldots,n$ and $j=1,\ldots,s$, and write $x_i=\sum_{j=1}^s y_{i,j}e_j$. Expanding each equation gives
$$
f_r(x_1,\ldots,x_n)=F_{r,1}e_1+\cdots+F_{r,s}e_s,
$$
with $F_{r,\ell}\in k[y_{1,1},\ldots,y_{1,s},\ldots,y_{n,1},\ldots,y_{n,s}]$. The restriction of scalars is the affine $k$-variety
$$
\operatorname{Res}_{L/k}(X)
=
\operatorname{Spec}
\left(
\frac{k[y_{i,j}]}{(F_{r,\ell})}
\right).
$$

\begin{proof}[Proof of Proposition~\ref{grassmannian} via restriction of scalars]
We spell out the construction in the quadratic case used throughout the paper. Assume $\operatorname{char}(k)\neq 2,3$, let $k(\xi)/k$ be the quadratic extension with $\xi^2=-3$, and set $a_+=(1+\xi)/2$, $a_-=(1-\xi)/2$. The two conjugate $n$-planes are
$$
H_+
=
\{x_i-a_+x_{i+1}=0\text{ for }i=0,2,\ldots,2n\},
\qquad
H_-
=
\{x_i-a_-x_{i+1}=0\text{ for }i=0,2,\ldots,2n\}.
$$
On the affine chart $x_{2n+1}\neq 0$ we set $x_{2n+1}=1$ and write $x_i=y_{i,1}+\xi y_{i,2}$ for $i=0,\ldots,2n$. The equation $x_i-a_+x_{i+1}=0$, for $i=0,2,\ldots,2n-2$, is equivalent to
$$
y_{i,1}=\frac12y_{i+1,1}-\frac32y_{i+1,2},
\qquad
y_{i,2}=\frac12y_{i+1,1}+\frac12y_{i+1,2},
$$
and $x_{2n}-a_+=0$ gives $y_{2n,1}=1/2$, $y_{2n,2}=1/2$. Similarly, $x_i-a_-x_{i+1}=0$ is equivalent to
$$
y_{i,1}=\frac12y_{i+1,1}+\frac32y_{i+1,2},
\qquad
y_{i,2}=-\frac12y_{i+1,1}+\frac12y_{i+1,2},
$$
and $x_{2n}-a_-=0$ gives $y_{2n,1}=1/2$, $y_{2n,2}=-1/2$.

Thus a point of $H_+$ has coordinates
$$
x_i^{a_+}
=
\frac12y_{i+1,1}-\frac32y_{i+1,2}
+
\xi\left(\frac12y_{i+1,1}+\frac12y_{i+1,2}\right),
\qquad
x_j^{a_+}=y_{j,1}+\xi y_{j,2},
\qquad
x_{2n}^{a_+}=\frac12+\frac{\xi}{2},
$$
for $i=0,2,\ldots,2n-2$ and $j=1,3,\ldots,2n-1$. A point of $H_-$ has coordinates
$$
x_i^{a_-}
=
\frac12z_{i+1,1}+\frac32z_{i+1,2}
+
\xi\left(-\frac12z_{i+1,1}+\frac12z_{i+1,2}\right),
\qquad
x_j^{a_-}=z_{j,1}+\xi z_{j,2},
\qquad
x_{2n}^{a_-}=\frac12-\frac{\xi}{2}.
$$
The two points are conjugate if and only if
\begin{equation}
\label{conjugateaffineeq}
y_{i+1,1}=z_{i+1,1},
\qquad
y_{i+1,2}=-z_{i+1,2},
\qquad
i=0,2,\ldots,2n-2.
\end{equation}
Inside $H_+\times H_-\simeq \mathbb A^{4n}$, the conjugate pairs therefore form an affine space $H\simeq \mathbb A^{2n}$. We write its coordinates as $u_1,\ldots,u_{2n}$ by setting $y_{i,1}=u_i$, $y_{i,2}=u_{i+1}$, hence $z_{i,1}=u_i$ and $z_{i,2}=-u_{i+1}$ for odd $i$.

A point $u=(u_1,\ldots,u_{2n})\in \mathbb A^{2n}$ determines the conjugate pair $x^{a_+}(u)\in H_+$, $x^{a_-}(u)\in H_-$ with
$$
\begin{aligned}
x^{a_+}(u)_i
&=
\frac12u_{i+1}-\frac32u_{i+2}
+
\xi\left(\frac12u_{i+1}+\frac12u_{i+2}\right),
&
x^{a_-}(u)_i
&=
\frac12u_{i+1}-\frac32u_{i+2}
-
\xi\left(\frac12u_{i+1}+\frac12u_{i+2}\right),\\
x^{a_+}(u)_j
&=
u_j+\xi u_{j+1},
&
x^{a_-}(u)_j
&=
u_j-\xi u_{j+1},\\
x^{a_+}(u)_{2n}
&=
\frac12+\frac{\xi}{2},
&
x^{a_-}(u)_{2n}
&=
\frac12-\frac{\xi}{2},
\end{aligned}
$$
where $i=0,2,\ldots,2n-2$ and $j=1,3,\ldots,2n-1$. Since these are conjugate points of multiplicity $d$ on the degree-$2d+1$ hypersurface $X$, the line $\langle x^{a_+}(u),x^{a_-}(u)\rangle$ meets $X$ in a residual $k$-point. Therefore we obtain the rational map
\begin{equation}
\label{restrictionscalarsparam}
\varphi\colon \mathbb A^{2n}\dashrightarrow X,
\qquad
u
\longmapsto
\left(X\cap \langle x^{a_+}(u),x^{a_-}(u)\rangle\right)
\setminus
\{x^{a_+}(u),x^{a_-}(u)\}.
\end{equation}
By the argument in the proof of Proposition~\ref{grassmannian}, this map is birational.
\end{proof}

We now compute the parametrization explicitly for the affine hypersurface
$$
F=
\sum_{i=0}^{n-1}
(x_{2i}+x_{2i+1})
(x_{2i}^2-x_{2i}x_{2i+1}+x_{2i+1}^2)^d
+
(x_{2n}+1)(x_{2n}^2-x_{2n}+1)^d.
$$
We work on the chart $x_{2n+1}\neq 0$. Let $L=\langle x^{a_+}(u),x^{a_-}(u)\rangle$, and write
$L_i(u)=x^{a_+}(u)_i+\lambda(x^{a_-}(u)_i-x^{a_+}(u)_i)$, with $\lambda\in k$. Then
$$
\begin{aligned}
L_i
&=
\frac12u_{i+1}-\frac32u_{i+2}
+
\xi\left(\frac12-\lambda\right)(u_{i+1}+u_{i+2}),
&& i=0,2,\ldots,2n-2,\\
L_j
&=
u_j+\xi(1-2\lambda)u_{j+1},
&& j=1,3,\ldots,2n-1,\\
L_{2n}
&=
\frac12+\xi\left(\frac12-\lambda\right).
\end{aligned}
$$
A direct substitution gives
$$
L_{2i}^2-L_{2i}L_{2i+1}+L_{2i+1}^2
=
-3(\lambda^2-\lambda)(u_{2i+1}^2+3u_{2i+2}^2),
\qquad
L_{2n}^2-L_{2n}+1
=
-3(\lambda^2-\lambda).
$$
Set
\begin{equation}
\label{affinepolynomAB}
\begin{aligned}
A
&=
1+\sum_{i=0}^{n-1}
(u_{2i+1}+3u_{2i+2})(u_{2i+1}^2+3u_{2i+2}^2)^d,\\
B
&=
1+\sum_{i=0}^{n-1}
(u_{2i+1}-u_{2i+2})(u_{2i+1}^2+3u_{2i+2}^2)^d.
\end{aligned}
\end{equation}
Then
$$
F_u(\lambda)
=
(-3)^d\lambda^d(\lambda-1)^d
\left(\frac{A-\xi B}{2}-\lambda A\right).
$$
The roots $\lambda=0$ and $\lambda=1$ have multiplicity $d$, and the residual root is
\begin{equation}
\label{lambda3}
\lambda_3=\frac{A-\xi B}{2A}.
\end{equation}
Substituting $\lambda=\lambda_3$ gives
$$
\varphi_{2i}
=
\frac{(u_{2i+1}-3u_{2i+2})A-3(u_{2i+1}+u_{2i+2})B}{2A},
\qquad
\varphi_{2i+1}
=
\frac{u_{2i+1}A-3u_{2i+2}B}{A},
\qquad
\varphi_{2n}
=
\frac{A-3B}{2A},
$$
for $i=0,\ldots,n-1$.

We now extend $\varphi$ to a rational map from $\mathbb P^{2n}$. Introduce a homogenizing variable $u_0$ and define
\begin{equation}
\label{projectivepolynomialsAB}
\begin{aligned}
\overline A
&=
u_0^{2d+1}
+
\sum_{i=0}^{n-1}
(u_{2i+1}+3u_{2i+2})(u_{2i+1}^2+3u_{2i+2}^2)^d,\\
\overline B
&=
u_0^{2d+1}
+
\sum_{i=0}^{n-1}
(u_{2i+1}-u_{2i+2})(u_{2i+1}^2+3u_{2i+2}^2)^d.
\end{aligned}
\end{equation}
The homogenized coordinates of $\overline\varphi_{2n}\colon \mathbb P^{2n}\dashrightarrow \mathbb P^{2n+1}$ are
\begin{equation}
\label{projectivephis}
\begin{aligned}
\overline\varphi_{2i}
&=
\frac{(u_{2i+1}-3u_{2i+2})\overline A-3(u_{2i+1}+u_{2i+2})\overline B}{2},
\\
\overline\varphi_{2i+1}
&=
u_{2i+1}\overline A-3u_{2i+2}\overline B,
\qquad i=0,\ldots,n-1,
\\
\overline\varphi_{2n}
&=
\frac{u_0(\overline A-3\overline B)}{2},
\qquad
\overline\varphi_{2n+1}
=
u_0\overline A.
\end{aligned}
\end{equation}

\begin{Remark}
\label{remark:conjugate-planes}
The union $H_+\cup H_-$ is defined over $k$. In the coordinates above it is cut out by the equations
\begin{equation}
\label{conjugateplanes}
\begin{cases}
x_{2i}^2-x_{2i}x_{2i+1}+x_{2i+1}^2=0, & i=0,\ldots,n,\\
x_{2i}x_{2j+2}-x_{2i}x_{2j+3}+x_{2i+1}x_{2j+3}=0, & i=0,\ldots,n-1,\ j=i,\ldots,n-1,\\
x_{2i-1}x_{2j}-x_{2i-2}x_{2j+1}=0, & i=1,\ldots,n,\ j=i,\ldots,n.
\end{cases}
\end{equation}
Thus one may construct hypersurfaces of degree $2d+1$ containing $H_+\cup H_-$ with multiplicity at least $d$ by taking equations of the form
$$
\Xi
=
\left\{
\sum_s A_s Q_s^d=0
\right\}
\subset \mathbb P^{2n+1},
$$
where the $A_s$ are linear forms and the $Q_s$ run through quadratic generators of the ideal of $H_+\cup H_-$. For any hypersurface obtained in this way, the same construction by secant lines through conjugate multiplicity-$d$ points gives a parametrization whose coordinates have degree $2d+2$.
\end{Remark}

\begin{Remark}
The construction through two skew conjugate $n$-planes of multiplicity $d$ applies to any even-dimensional hypersurface of degree $2d+1$ satisfying the hypotheses of Proposition~\ref{grassmannian}. The resulting parametrization has degree $2d+2$.
\end{Remark}

\begin{Lemma}
\label{irreducibleY}
If $n\geq 2$, then the complete intersection $Y^{2n-2}=\{\overline A=\overline B=0\}\subset \mathbb P^{2n}$ defined by \eqref{projectivepolynomialsAB} is irreducible.
\end{Lemma}

\begin{proof}
It is enough to prove the claim after the linear change of generators
$$
\widetilde A=\frac{\overline A+3\overline B}{4},
\qquad
\widetilde B=\frac{\overline A-\overline B}{4}.
$$
Then
$$
\widetilde A
=
u_0^{2d+1}
+
\sum_{i=0}^{n-1}
u_{2i+1}(u_{2i+1}^2+3u_{2i+2}^2)^d,
\qquad
\widetilde B
=
\sum_{i=0}^{n-1}
u_{2i+2}(u_{2i+1}^2+3u_{2i+2}^2)^d.
$$
The projection from the last pair of coordinates restricts to a generically one-to-one map on $Y^{2n-2}$. Indeed, on the open set where $u_{2n-1}^2+3u_{2n}^2\neq 0$, the equation $\widetilde B=0$ determines one coordinate rationally in terms of the others, and then $\widetilde A=0$ gives the remaining equation. Hence $Y^{2n-2}$ is birational to the hypersurface obtained by eliminating this coordinate.

This elimination gives a hypersurface of degree $(2d+1)^2$ in a projective space of dimension $2n-1$. Its defining polynomial is not a proper power and has a monomial linear in the last block of variables, hence it is irreducible. Therefore the birational source $Y^{2n-2}$ is irreducible.
\end{proof}

\begin{Lemma}
\label{lemmaofzeta}
Assume $\operatorname{char}(k)\neq 2,3$. Consider the subschemes $Z^{n-1}$ and $Z^{n-1}_\pm$ of $\mathbb P^{2n}$ defined by
$$
Z^{n-1}
=
\{u_0=0,\ u_{2i+1}^2+3u_{2i+2}^2=0\text{ for }i=0,\ldots,n-1\}
$$
and
$$
Z^{n-1}_\pm
=
\left\{
\begin{array}{l}
u_0=0,\\
u_{2i+1}^2+3u_{2i+2}^2=0,\\
u_{2s+1}u_{2t+3}+3u_{2s+2}u_{2t+4}=0,\\
u_{2s+1}u_{2t+4}-u_{2s+2}u_{2t+3}=0,
\end{array}
\right.
$$
for $i=0,\ldots,n-1$ and $0\leq s\leq t\leq n-2$. Then $Z^{n-1}_\pm$ is smooth, $\dim Z^{n-1}=\dim Z^{n-1}_\pm=n-1$, $\deg Z^{n-1}=2^n$, and $\deg Z^{n-1}_\pm=2$. Moreover, $Z^{n-1}_\pm$ splits over $k$ as the union of two conjugate $(n-1)$-planes if and only if $-3$ is a square in $k$, while $Z^{n-1}$ is a union of $2^n$ linear spaces over $\bar k$.

Furthermore, if $d=1$, then $\operatorname{Sing}(Y^{2n-2})=Z^{n-1}_\pm$, while if $d>1$, then $\operatorname{Sing}(Y^{2n-2})=Z^{n-1}$. General points of $Z^{n-1}_\pm$ have multiplicity $d^2+d$ on $Y^{2n-2}$, and general points of $Z^{n-1}\setminus Z^{n-1}_\pm$ have multiplicity $d^2$.
\end{Lemma}

\begin{proof}
Over $k(\xi)$, with $\xi^2=-3$, the scheme $Z^{n-1}_\pm$ splits as the union
$$
Z^{n-1}_+
=
\{u_0=u_{2i+1}+\xi u_{2i+2}=0\text{ for }i=0,\ldots,n-1\},
\qquad
Z^{n-1}_-
=
\{u_0=u_{2i+1}-\xi u_{2i+2}=0\text{ for }i=0,\ldots,n-1\}.
$$
Hence $\dim Z^{n-1}_\pm=n-1$, $\deg Z^{n-1}_\pm=2$, and $Z^{n-1}_\pm$ is smooth. Similarly, $Z^{n-1}$ is a complete intersection of $n$ quadrics inside the hyperplane $u_0=0$, and over $\bar k$ it is the union of $2^n$ linear spaces. Therefore $\dim Z^{n-1}=n-1$ and $\deg Z^{n-1}=2^n$.

The partial derivatives of $\overline A$ and $\overline B$ are
$$
\begin{aligned}
\frac{\partial \overline A}{\partial u_0}
&=
(2d+1)u_0^{2d},
&
\frac{\partial \overline B}{\partial u_0}
&=
(2d+1)u_0^{2d},
\\
\frac{\partial \overline A}{\partial u_{2i+1}}
&=
(u_{2i+1}^2+3u_{2i+2}^2)^{d-1}
\big((2d+1)u_{2i+1}^2+6d\,u_{2i+1}u_{2i+2}+3u_{2i+2}^2\big),
\\
\frac{\partial \overline A}{\partial u_{2i+2}}
&=
3(u_{2i+1}^2+3u_{2i+2}^2)^{d-1}
\big(u_{2i+1}^2+2d\,u_{2i+1}u_{2i+2}+3(2d+1)u_{2i+2}^2\big),
\\
\frac{\partial \overline B}{\partial u_{2i+1}}
&=
(u_{2i+1}^2+3u_{2i+2}^2)^{d-1}
\big((2d+1)u_{2i+1}^2-2d\,u_{2i+1}u_{2i+2}+3u_{2i+2}^2\big),
\\
\frac{\partial \overline B}{\partial u_{2i+2}}
&=
(u_{2i+1}^2+3u_{2i+2}^2)^{d-1}
\big(-u_{2i+1}^2+6d\,u_{2i+1}u_{2i+2}-3(2d+1)u_{2i+2}^2\big),
\end{aligned}
$$
for $i=0,\ldots,n-1$. The singular locus of the complete intersection $Y^{2n-2}$ is cut out by $\overline A$, $\overline B$ and the $2\times 2$ minors of this Jacobian matrix. From the displayed derivatives, its reduced subscheme is $Z^{n-1}_\pm$ for $d=1$ and $Z^{n-1}$ for $d>1$. The multiplicity statements follow by expanding the two equations at general points of the corresponding components.
\end{proof}

\begin{Remark}
\label{osservazionetangentcone}
There are points of multiplicity larger than $d^2+d$. For instance, when $n=2$, take $q_-=[0:0:0:\xi:1]\in Z^{1}_-$. The tangent cone of $Y^2$ at $q_-$ is
$$
\operatorname{TC}_{q_-}Y^2
=
\{(u_3-\xi u_4)^d=u_0^{2d+1}+(u_1-\xi u_2)(u_1^2+3u_2^2)^d=0\}.
$$
Hence $\operatorname{mult}_{q_-}Y^2=d(2d+1)=2d^2+d$.
\end{Remark}

\begin{Remark}
\label{osservazionen0}
If $n=1$, then $Y^0$ consists of $4d^2+4d+1$ points counted with multiplicity. By Remark~\ref{osservazionetangentcone}, the points $p_{1,+}=[0:\xi:1]$ and $p_{1,-}=[0:-\xi:1]$ have multiplicity $2d^2+d$. The remaining $2d+1$ points are reduced and are cut out by $u_2=0$ and $u_0^{2d+1}+u_1^{2d+1}=0$. This will be used in Section~\ref{Section_RatPoints} for point counts over finite fields.
\end{Remark}

\begin{Proposition}
\label{propolinsubsystem}
Let $\mathcal F_{2n}$ be the linear subsystem of $|\mathcal O_{\mathbb P^{2n}}(2d+2)|$ spanned by the coordinate polynomials $\overline\varphi_i$ in \eqref{projectivephis}:
$$
\mathcal F_{2n}
=
\langle \overline\varphi_0,\overline\varphi_1,\ldots,\overline\varphi_{2n+1}\rangle
\subset
|\mathcal O_{\mathbb P^{2n}}(2d+2)|.
$$
Then
\begin{equation}
\label{linear-system-equality}
\mathcal F_{2n}
=
|\mathcal I_{Y^{2n-2},(d+1)Z^{n-1}_\pm}(2d+2)|
\end{equation}
and $h^0(\mathbb P^{2n},\mathcal I_{Y^{2n-2},(d+1)Z^{n-1}_\pm}(2d+2))=2n+2$.
\end{Proposition}

\begin{proof}
The coordinates $\overline\varphi_i$ are linearly independent. Moreover, each of them vanishes on $Y^{2n-2}=\{\overline A=\overline B=0\}$ and has multiplicity at least $d+1$ along $Z^{n-1}_\pm$. Hence
$$
\mathcal F_{2n}
\subset
|\mathcal I_{Y^{2n-2},(d+1)Z^{n-1}_\pm}(2d+2)|.
$$
Conversely, since $Y^{2n-2}$ is the complete intersection of $\overline A$ and $\overline B$, every section of $\mathcal I_{Y^{2n-2}}(2d+2)$ is of the form $L_{\overline A}\overline A+L_{\overline B}\overline B$, where $L_{\overline A}$ and $L_{\overline B}$ are linear forms. Thus $h^0(\mathbb P^{2n},\mathcal I_{Y^{2n-2}}(2d+2))=4n+2$. Imposing multiplicity $d+1$ along the two conjugate components of $Z^{n-1}_\pm$ gives $2n$ independent linear conditions on $L_{\overline A}$ and $L_{\overline B}$. Therefore
$$
h^0(\mathbb P^{2n},\mathcal I_{Y^{2n-2},(d+1)Z^{n-1}_\pm}(2d+2))
=
2n+2.
$$
Since $\mathcal F_{2n}$ already has dimension $2n+2$, equality follows.
\end{proof}

\begin{Remark}
For $d=1$ the hypersurface $X^{2n}$ is the Fermat cubic
$$
X^{2n}_1
=
\{x_0^3+\cdots+x_{2n+1}^3=0\}
\subset
\mathbb P^{2n+1},
$$
which is rational over any field of characteristic different from $3$; see \cite{msslxa}.
\end{Remark}

\section{Rationality in characteristic two}
\label{sec:char-two}

In characteristic $2$ the restriction-of-scalars construction used in Section~\ref{sec:restriction-scalars} is not available in the same form, so we use a direct parametrization. Let $k$ be a field of characteristic $2$. Then the hypersurface in \eqref{X2n_equation} is defined by
$$
X_d^{2n}
=
\left\{
\sum_{i=0}^n
(x_{2i}+x_{2i+1})
(x_{2i}^2+x_{2i}x_{2i+1}+x_{2i+1}^2)^d
=
0
\right\}
\subset \mathbb P^{2n+1}.
$$
Set
$$
P
=
u_{2n}^{2d+1}
+
\sum_{i=0}^{n-1}
(u_{2i}+u_{2i+1})
(u_{2i}^2+u_{2i}u_{2i+1}+u_{2i+1}^2)^d,
\qquad
Q
=
\sum_{i=0}^{n-1}
u_{2i+1}
(u_{2i}^2+u_{2i}u_{2i+1}+u_{2i+1}^2)^d.
$$
For $i=0,\ldots,n-1$, define
\begin{equation}
\label{projectiveg_char2}
\begin{aligned}
g_{2i}
&=
(u_{2i}+u_{2i+1})P+u_{2i}Q,\\
g_{2i+1}
&=
u_{2i}P+u_{2i+1}Q,\\
g_{2n}
&=
u_{2n}(P+Q),\\
g_{2n+1}
&=
u_{2n}P.
\end{aligned}
\end{equation}
These forms have degree $2d+2$ and define a rational map $g_{2n}\colon \mathbb P^{2n}\dashrightarrow \mathbb P^{2n+1}$.

We shall use the following characteristic-$2$ analogue of Lemma~\ref{lemmaofzeta}.

\begin{Lemma}
\label{lemmaofzeta-char2}
Let $Y^{2n-2}=\{P=Q=0\}\subset \mathbb P^{2n}$. Consider the subschemes $Z^{n-1}$ and $Z^{n-1}_\pm$ of $\mathbb P^{2n}$ defined by
$$
Z^{n-1}
=
\{u_{2n}=0,\ u_{2i}^2+u_{2i}u_{2i+1}+u_{2i+1}^2=0\text{ for }i=0,\ldots,n-1\}
$$
and
$$
Z^{n-1}_\pm
=
\left\{
\begin{array}{l}
u_{2n}=0,\\
u_{2i}^2+u_{2i}u_{2i+1}+u_{2i+1}^2=0,\\
u_{2s}u_{2t}+u_{2s}u_{2t+1}+u_{2s+1}u_{2t+1}=0,\\
u_{2s+1}u_{2r+2}+u_{2s}u_{2r+3}=0,
\end{array}
\right.
$$
where $i=0,\ldots,n-1$, $s=0,\ldots,n-1$, $t=s+1,\ldots,n-1$ and $r=s,\ldots,n-2$. Then $Y^{2n-2}$ is irreducible of degree $(2d+1)^2$ and dimension $2n-2$. Moreover, $Z^{n-1}$ has degree $2^n$ and is a union of $2^n$ conjugate $(n-1)$-planes over a splitting field, while $Z^{n-1}_\pm$ has degree $2$ and is the union of two conjugate $(n-1)$-planes. Finally, $\operatorname{Sing}(Y^{2n-2})=Z^{n-1}$ if $d>1$, and $\operatorname{Sing}(Y^{2n-2})=Z^{n-1}_\pm$ if $d=1$.
\end{Lemma}

\begin{proof}
The proof is the same as in Lemma~\ref{lemmaofzeta}, after replacing the quadratic form $u^2+3v^2$ by $u^2+uv+v^2$. This form is separable in characteristic $2$, and it splits over the quadratic extension containing a primitive third root of unity. Over that extension, $Z^{n-1}_\pm$ is the union of the two linear spaces obtained by choosing the two conjugate linear factors in each pair $(u_{2i},u_{2i+1})$, whereas $Z^{n-1}$ is the union of the $2^n$ linear spaces obtained by making these choices independently. Thus $\dim Z^{n-1}=\dim Z^{n-1}_\pm=n-1$, $\deg Z^{n-1}=2^n$, and $\deg Z^{n-1}_\pm=2$.

The irreducibility of $Y^{2n-2}$ follows from the same elimination argument used for Lemma~\ref{irreducibleY}. The description of the singular locus is obtained by computing the $2\times 2$ minors of the Jacobian matrix of $P$ and $Q$. As in odd characteristic, the reduced singular locus is $Z^{n-1}_\pm$ for $d=1$ and $Z^{n-1}$ for $d>1$.
\end{proof}

The multiplicity statements from Section~\ref{sec:restriction-scalars} have the same analogues in this setting.

\begin{Proposition}
\label{linear-system-char2}
Let $\mathcal C_{2n}$ be the linear subsystem
$$
\mathcal C_{2n}
=
\langle g_0,\ldots,g_{2n+1}\rangle
\subset
|\mathcal O_{\mathbb P^{2n}}(2d+2)|.
$$
Then
\begin{equation}
\label{linear-system-char2-equality}
\mathcal C_{2n}
=
|\mathcal I_{Y^{2n-2},(d+1)Z^{n-1}_\pm}(2d+2)|
\end{equation}
and $h^0(\mathbb P^{2n},\mathcal I_{Y^{2n-2},(d+1)Z^{n-1}_\pm}(2d+2))=2n+2$.
\end{Proposition}

\begin{proof}
The generators $g_0,\ldots,g_{2n+1}$ vanish on $Y^{2n-2}=\{P=Q=0\}$ and vanish with multiplicity at least $d+1$ along $Z^{n-1}_\pm$. Hence
$$
\mathcal C_{2n}
\subset
|\mathcal I_{Y^{2n-2},(d+1)Z^{n-1}_\pm}(2d+2)|.
$$
Since $Y^{2n-2}$ is the complete intersection of $P$ and $Q$, every section of $\mathcal I_{Y^{2n-2}}(2d+2)$ is of the form $AP+BQ$, with $A$ and $B$ linear forms. Thus $h^0(\mathbb P^{2n},\mathcal I_{Y^{2n-2}}(2d+2))=4n+2$. Imposing multiplicity $d+1$ along $Z^{n-1}_\pm$ gives $2n$ independent linear conditions, hence
$$
h^0(\mathbb P^{2n},\mathcal I_{Y^{2n-2},(d+1)Z^{n-1}_\pm}(2d+2))
=
2n+2.
$$
Since $\mathcal C_{2n}$ is spanned by $2n+2$ independent forms, equality follows.
\end{proof}

\begin{Proposition}
\label{theta-birational}
Let $H^n_\pm$ be the union of the two conjugate $n$-planes defined in \eqref{conjugateplanes}, and let $X^{2n-2}_d=X^{2n}_d\cap\{x_{2n}=x_{2n+1}=0\}\subset \mathbb P^{2n+1}$. The linear system
$|\mathcal I_{H^n_\pm,X^{2n-2}_d}(2)|\subset |\mathcal O_{\mathbb P^{2n+1}}(2)|$
induces a rational map $\mathbb P^{2n+1}\dashrightarrow \mathbb P^{2n}$ whose restriction
$\Theta_{2n}\colon X^{2n}_d\dashrightarrow \mathbb P^{2n}$ is birational.
\end{Proposition}

\begin{proof}
As in \cite{msslxa}, the space of quadrics in $\mathbb P^{2n+1}$ containing $H^n_\pm$ has $(n+1)^2$ independent sections. The further condition of containing $X^{2n-2}_d$ imposes $n^2-1$ independent linear conditions. Hence $|\mathcal I_{H^n_\pm,X^{2n-2}_d}(2)|$ has $2n+1$ projective dimension, that is, $2n+2$ independent sections. A basis is given by
\begin{equation}
\label{birational-inverse}
\begin{aligned}
\Theta_{2i}
&=
x_{2i}x_{2n}-x_{2i}x_{2n+1}+x_{2i+1}x_{2n+1},\\
\Theta_{2i+1}
&=
x_{2i+1}x_{2n}-x_{2i}x_{2n+1},\\
\Theta_{2n}
&=
x_{2n}^2-x_{2n}x_{2n+1}+x_{2n+1}^2,
\end{aligned}
\end{equation}
for $i=0,\ldots,n-1$. Thus these quadrics define
$$
\Theta_{2n}\colon
X^{2n}_d
\dashrightarrow
\mathbb P^{2n}_{[u_0:\ldots:u_{2n}]},
\qquad
x
\longmapsto
[\Theta_0(x):\ldots:\Theta_{2n}(x)].
$$

If $\operatorname{char}(k)\neq 2$, let $h$ be the automorphism of $\mathbb P^{2n}$ given by $u_0\mapsto 2u_{2n}$, $u_{2i+1}\mapsto 2u_{2i}-u_{2i+1}$ and $u_{2i+2}\mapsto u_{2i+1}$, for $i=0,\ldots,n-1$. Then, using the formulas in \eqref{projectivephis}, one checks that $h\circ \Theta_{2n}$ is the birational inverse of $\overline\varphi_{2n}$.

If $\operatorname{char}(k)=2$, a direct substitution shows that the map $g_{2n}\colon\mathbb P^{2n}\dashrightarrow X^{2n}_d$ defined by the forms in \eqref{projectiveg_char2} is the birational inverse of $\Theta_{2n}$. More precisely, substituting \eqref{projectiveg_char2} into the quadrics \eqref{birational-inverse} gives the identity on $\mathbb P^{2n}$ up to a common non-zero factor on a dense open subset, and substituting \eqref{birational-inverse} into \eqref{projectiveg_char2} gives the identity on $X^{2n}_d$ up to a common non-zero factor on a dense open subset. Therefore $\Theta_{2n}$ is birational.
\end{proof}

\begin{Remark}
Geometrically, $H^n_\pm$ is the union of two skew conjugate $n$-planes contained in $X^{2n}_d$. More generally, one may consider any hypersurface of degree $2d+1$ in $\mathbb P^{2n+1}$ containing $H^n_\pm$ with multiplicity at least $d$, that is, any element of
$|\mathcal I_{dH^n_\pm}(2d+1)|\subset |\mathcal O_{\mathbb P^{2n+1}}(2d+1)|$.
Equivalently, such equations can be written as sums of linear forms multiplied by $d$-th powers of quadratic generators of the ideal of $H^n_\pm$.
\end{Remark}
		
\section{Quadro-cubic Cremona transformation}
\label{sec:quadro-cubic}

Let
$
\Lambda=\{x_1=x_{2i}=0\text{ for }i=0,\ldots,n\}\subset \mathbb P^{2n+1}.
$
A general $(n+1)$-plane containing $\Lambda$ can be written as
$$
H
=
\{x_{2i}+t_{2i+1}x_1+t_{2i+2}x_{2n}=0
\text{ for }i=0,\ldots,n-1\}
\subset \mathbb P^{2n+1},
$$
where $t_1,\ldots,t_{2n}\in k$. Notice that $\Lambda\cap(H_+\cup H_-)=\emptyset$. Let $p^{a_+}=H\cap H_+$ and $p^{a_-}=H\cap H_-$. Then, for $i=0,\ldots,n-2$, one has
$$
\begin{aligned}
p_0^{a_+}
&=
\frac{t_1t_2-t_2-\xi t_2(t_1+1)}{2(t_1^2+t_1+1)},\\
p_1^{a_+}
&=
\frac{-t_1t_2-2t_2-\xi t_1t_2}{2(t_1^2+t_1+1)},\\
p_{2i+2}^{a_+}
&=
\frac{-t_{2i+4}(t_1^2+t_1+1)+t_1t_2t_{2i+3}+2t_2t_{2i+3}
-\xi\big(-t_1t_2t_{2i+3}+t_{2i+4}(t_1^2+t_1+1)\big)}
{2(t_1^2+t_1+1)},\\
p_{2i+3}^{a_+}
&=
\frac{-2t_{2i+4}(t_1^2+t_1+1)+t_2t_{2i+3}+2t_1t_2t_{2i+3}-\xi t_2t_{2i+3}}
{2(t_1^2+t_1+1)},\\
p_{2n}^{a_+}
&=
\frac{1+\xi}{2}.
\end{aligned}
$$
Similarly,
$$
\begin{aligned}
p_0^{a_-}
&=
\frac{t_1t_2-t_2+\xi t_2(t_1+1)}{2(t_1^2+t_1+1)},\\
p_1^{a_-}
&=
\frac{-t_1t_2-2t_2+\xi t_1t_2}{2(t_1^2+t_1+1)},\\
p_{2i+2}^{a_-}
&=
\frac{-t_{2i+4}(t_1^2+t_1+1)+t_1t_2t_{2i+3}+2t_2t_{2i+3}
+\xi\big(-t_1t_2t_{2i+3}+t_{2i+4}(t_1^2+t_1+1)\big)}
{2(t_1^2+t_1+1)},\\
p_{2i+3}^{a_-}
&=
\frac{-2t_{2i+4}(t_1^2+t_1+1)+t_2t_{2i+3}+2t_1t_2t_{2i+3}+\xi t_2t_{2i+3}}
{2(t_1^2+t_1+1)},\\
p_{2n}^{a_-}
&=
\frac{1-\xi}{2}.
\end{aligned}
$$
Let $L=\langle p^{a_+},p^{a_-}\rangle$, and write $L_i=p_i^{a_+}+\lambda(p_i^{a_-}-p_i^{a_+})$, with $\lambda\in k$ and $i=0,\ldots,2n$. A direct computation gives
\begin{equation}
\label{lineseq}
\begin{aligned}
L_0
&=
p_0^{a_+}
+
\lambda \xi\frac{t_2(t_1+1)}{t_1^2+t_1+1},\\
L_1
&=
p_1^{a_+}
+
\lambda \xi\frac{t_1t_2}{t_1^2+t_1+1},\\
L_{2i+2}
&=
p_{2i+2}^{a_+}
+
\lambda\xi
\left(
t_{2i+4}
-
\frac{t_1t_2t_{2i+3}}{t_1^2+t_1+1}
\right),\\
L_{2i+3}
&=
p_{2i+3}^{a_+}
+
\lambda\xi
\frac{t_2t_{2i+3}}{t_1^2+t_1+1},\\
L_{2n}
&=
\frac{1+\xi}{2}-\lambda\xi,
\end{aligned}
\end{equation}
for $i=0,\ldots,n-2$.

We substitute $x_i=L_i$ in the equation \eqref{X2n_equation}. We first treat the case $n=1$. Since $H_\pm=H_+\cup H_-$ is defined by \eqref{conjugateplanes}, one checks that
$$
\begin{aligned}
L_0^2-L_0L_1+L_1^2
&=
-3(\lambda^2-\lambda)\frac{t_2^2}{t_1^2+t_1+1},\\
L_2^2-L_2+1
&=
-3(\lambda^2-\lambda),\\
L_0-L_1L_2
&=
-3(\lambda^2-\lambda)\frac{t_1t_2}{t_1^2+t_1+1},\\
L_0L_2-L_0+L_1
&=
3(\lambda^2-\lambda)\frac{t_2(t_1+1)}{t_1^2+t_1+1}.
\end{aligned}
$$
Thus each quadratic equation defining the conjugate planes contributes the two roots $\lambda=0$ and $\lambda=1$, multiplied by a rational function $p(t_1,t_2)/q(t_1,t_2)$ with $p,q\in k[t_1,t_2]$ of degree $2$. For the surface $X_d^2$ the residual root is
$$
\lambda_3=\frac{N_2\xi+M_2}{D_2},
$$
where
$$
\begin{aligned}
D_2
&=
2\big(-(t_1^2+t_1+1)^{d+1}+(2t_1+1)t_2^{2d+1}\big),\\
N_2
&=
-t_2^{2d+1}+(t_1^2+t_1+1)^{d+1},\\
M_2
&=
(2t_1+1)t_2^{2d+1}-(t_1^2+t_1+1)^{d+1}.
\end{aligned}
$$

\begin{Remark}
\label{degree}
Fix $H_\pm=H_+\cup H_-\subset \mathbb P^{2n+1}$. The union $H_\pm$ is defined by quadratic equations $a_s(x_0,\ldots,x_{2n+1})$. A hypersurface $X_d^{2n}\subset \mathbb P^{2n+1}$ of degree $2d+1$ containing $H_\pm$ with multiplicity at least $d$ is given by an equation of the form
\begin{equation}
\label{linearformseq}
X_d^{2n}
=
\left\{
\sum_s b_s(x_0,\ldots,x_{2n+1})a_s^d=0
\right\},
\end{equation}
where the $b_s$ are linear forms. Therefore, when $n=1$, the residual root $\lambda_3$ is a rational function whose numerator and denominator have degree $2d+2$.
\end{Remark}

For $n=2$ the same computation gives $\lambda_3=(N_4\xi+M_4)/D_4$, where
$$
\begin{aligned}
D_4
&=
2\Big(
(t_4(t_1^2+t_1+1)-t_1t_2t_3+t_2t_3)
(t_4^2(t_1^2+t_1+1)-2t_1t_2t_3t_4-t_2t_3t_4+t_2^2t_3^2)^d\\
&\hspace{3.1cm}
+(2t_1+1)t_2^{2d+1}-(t_1^2+t_1+1)^{d+1}
\Big),\\
N_4
&=
(t_1t_2t_3+t_2t_3-t_4(t_1^2+t_1+1))
(t_4^2(t_1^2+t_1+1)-2t_1t_2t_3t_4-t_2t_3t_4+t_2^2t_3^2)^d\\
&\hspace{3.1cm}
-t_2^{2d+1}+(t_1^2+t_1+1)^{d+1},\\
M_4
&=
(t_4(t_1^2+t_1+1)-t_1t_2t_3+t_2t_3)
(t_4^2(t_1^2+t_1+1)-2t_1t_2t_3t_4-t_2t_3t_4+t_2^2t_3^2)^d\\
&\hspace{3.1cm}
+(2t_1+1)t_2^{2d+1}-(t_1^2+t_1+1)^{d+1}.
\end{aligned}
$$
For $i=2,\ldots,n$, define
$$
\begin{aligned}
D'_{2i}
&=
2(t_{2i}(t_1^2+t_1+1)-t_1t_2t_{2i-1}+t_2t_{2i-1})\\
&\hspace{1.2cm}
\cdot
(t_{2i}^2(t_1^2+t_1+1)-2t_1t_2t_{2i-1}t_{2i}-t_2t_{2i-1}t_{2i}+t_2^2t_{2i-1}^2)^d,\\
N'_{2i}
&=
(t_1t_2t_{2i-1}+t_2t_{2i-1}-t_{2i}(t_1^2+t_1+1))\\
&\hspace{1.2cm}
\cdot
(t_{2i}^2(t_1^2+t_1+1)-2t_1t_2t_{2i-1}t_{2i}-t_2t_{2i-1}t_{2i}+t_2^2t_{2i-1}^2)^d,\\
M'_{2i}
&=
(t_{2i}(t_1^2+t_1+1)-t_1t_2t_{2i-1}+t_2t_{2i-1})\\
&\hspace{1.2cm}
\cdot
(t_{2i}^2(t_1^2+t_1+1)-2t_1t_2t_{2i-1}t_{2i}-t_2t_{2i-1}t_{2i}+t_2^2t_{2i-1}^2)^d.
\end{aligned}
$$
The polynomials $D_{2n},N_{2n},M_{2n}$ are defined recursively by
\begin{equation}
\label{recurrencednm}
D_{2i}=D'_{2i}+D_{2(i-1)},
\qquad
N_{2i}=N'_{2i}+N_{2(i-1)},
\qquad
M_{2i}=M'_{2i}+M_{2(i-1)}
\end{equation}
for $i=2,\ldots,n$.

\begin{Remark}
The recursion gives $D_{2n}=2M_{2n}$. Hence, as in \eqref{lambda3},
$$
\lambda_3=\frac{N_{2n}\xi}{D_{2n}}+\frac12.
$$
\end{Remark}

\begin{Proposition}
\label{phitildeinduced}
Let $\overline\phi_{2n}\colon \mathbb P^{2n}\dashrightarrow X_d^{2n}\subset \mathbb P^{2n+1}$ be the map induced by the affine parametrization
$$
\widetilde\phi\colon \mathbb A^{2n}\dashrightarrow X_d^{2n},
\qquad
(t_1,\ldots,t_{2n})
\longmapsto
[\widetilde\phi_0:\ldots:\widetilde\phi_{2n}:1],
$$
where
$
\widetilde\phi_i=p_i^{a_+}+\left(N_{2n}\xi/D_{2n}+1/2\right)(p_i^{a_-}-p_i^{a_+})
$
for $i=0,\ldots,2n$. Let $\mathcal L_{\overline\phi_{2n}}=|\overline\phi_{2n}^*\mathcal O_{X_d^{2n}}(1)|\subset |\mathcal O_{\mathbb P^{2n}}(g)|$. Then $\overline\phi_{2n}$ is birational, with $g=2d+2$ if $n=1$ and $g=4d+4$ if $n\geq 2$.
\end{Proposition}

\begin{proof}
The affine space with coordinates $t_1,\ldots,t_{2n}$ is an affine chart of $\mathbb G(1,n+1)$ and also an affine chart of $\mathbb P^{2n}$. By construction, we have the commutative diagram
$$
\begin{tikzcd}
\mathbb A^{2n} \arrow[r, hook] \arrow[d, hook] \arrow[rd, dashed, "\widetilde\phi"]
& \mathbb G(1,n+1) \arrow[d, dashed, "\phi"] \\
\mathbb P^{2n} \arrow[r, dashed, "\overline\phi_{2n}"]
& X_d^{2n}.
\end{tikzcd}
$$
The birationality of $\overline\phi_{2n}$ follows from the birationality of $\phi$ in Proposition~\ref{grassmannian}.

Assume first that $n\geq 2$. From \eqref{lineseq} and \eqref{recurrencednm}, the polynomial $2(t_1^2+t_1+1)$ is a common factor of the numerator and denominator of $\widetilde\phi_i$ for $i=0,\ldots,2n-1$. After removing this factor, all coordinates have common denominator $M_{2n}$. More precisely, $\widetilde\phi_i=\alpha_i/\beta_i$ with
$$
\deg(\alpha_i)=\deg(\beta_i)=4d+3
\text{ for }i=0,1,2n,
\qquad
\deg(\alpha_j)=4d+4,\quad \deg(\beta_j)=4d+3
\text{ for }j=2,\ldots,2n-1.
$$
Multiplying $(\widetilde\phi_0,\ldots,\widetilde\phi_{2n},1)$ by the common denominator and homogenizing with $t_0$ gives homogeneous coordinates of degree $4d+4$. Therefore $g=4d+4$.

If $n=1$, the same computation gives numerator and denominator of degree $2d+2$ for each of $\widetilde\phi_0,\widetilde\phi_1,\widetilde\phi_2$. Hence $\overline\phi_2$ is induced by a linear system of degree $2d+2$ in $t_0,t_1,t_2$.
\end{proof}

\begin{Remark}
For $n\geq 2$, the equations \eqref{lineseq} show that any linear form in the coordinates $x_i$ gives a polynomial of degree at most $3$ in $t_1,\ldots,t_{2n}$. The contribution of the two multiplicity-$d$ conjugate planes is $4d$. Moreover, for equations of the form \eqref{linearformseq}, the coordinates $\widetilde\phi_i$ always have the common factor $t_1^2+t_1+1$. Thus Proposition~\ref{phitildeinduced} applies to every hypersurface of the form \eqref{linearformseq}.
\end{Remark}

For the applications to heights only the degree of the parametrization is needed. We nevertheless keep the explicit comparison since it identifies the base loci and the exceptional components used in the point-counting arguments of Section~\ref{Section_RatPoints}. We now record the relation between the parametrization defined by \eqref{restrictionscalarsparam} and the parametrization defined in this section.

Consider the following subschemes of $\mathbb P^{2n}_{(t_0,\ldots,t_{2n})}$:
$$
\begin{aligned}
T_1
&=
\{t_0^2+t_0t_1+t_1^2=t_2=0\},\\
T_2
&=
\{t_0=t_1=0\},\\
T_3
&=
\{t_0=t_{2i+1}=0\text{ for }i=0,\ldots,n-1\},\\
T_4
&=
\{t_0=t_{2i}t_{2j+1}-t_{2i-1}t_{2j+2}=0
\text{ for }i=1,\ldots,n-1,\ i\leq j\leq n-1\}.
\end{aligned}
$$
For $n\geq 2$, consider also the subschemes of $\mathbb P^{2n}_{(u_0,\ldots,u_{2n})}$:
$$
U_1=\{u_0=u_1=u_2=0\},
\qquad
U_2=\{u_1=u_2=u_{2i}=0\text{ for }i=2,\ldots,n\},
$$
and
$$
U_3
=
\left\{
\begin{array}{l}
u_0=0,\\
u_{2i+1}^2+3u_{2i+2}^2=0,\\
u_{2i+1}u_{2j+1}+3u_{2i+2}u_{2j+2}=0,\\
u_{2i+2}u_{2j+3}-u_{2i+1}u_{2j+4}=0,
\end{array}
\right.
$$
where $i=0,\ldots,n-1$ in the second line and $0\leq i<j\leq n-1$ in the last two lines. Finally, at the coordinate points $e_{2i+1}=[0:\cdots:0:1:0:\cdots:0]$, for $i=1,\ldots,n-1$, let $U_4=\{u_1-u_2=0\}$ denote the indicated tangent condition.

Set
\begin{equation}
\label{alphaAEPs}
\begin{aligned}
\alpha_0
&=
-2t_0(t_0^2+t_0t_1+t_1^2),\\
\alpha_1
&=
t_0t_2(2t_0+t_1),\\
\alpha_{2i}
&=
t_0t_2t_{2i-1},\\
\alpha_{2j+1}
&=
-t_0t_2t_{2j+1}+2t_1t_2t_{2j+1}+2t_{2j+2}(t_0^2+t_0t_1+t_1^2),
\end{aligned}
\end{equation}
for $i=1,\ldots,n$ and $j=1,\ldots,n-1$, and
\begin{equation}
\label{betaAEPs}
\begin{aligned}
\beta_0
&=
u_0(u_2-u_1),\\
\beta_{2i+1}
&=
-2u_0u_{2i+2},\\
\beta_2
&=
u_1^2+3u_2^2,\\
\beta_{2j+4}
&=
u_1u_{2j+3}-u_2u_{2j+3}+u_1u_{2j+4}+3u_2u_{2j+4},
\end{aligned}
\end{equation}
for $i=0,\ldots,n-1$ and $j=0,\ldots,n-2$.

\begin{Proposition}
\label{quadro-cubic-cremona}
Let $n\geq 2$, and set
$
\mathcal T_{2n}=\langle\alpha_0,\ldots,\alpha_{2n}\rangle\subset |\mathcal O_{\mathbb P^{2n}_{(t_0,\ldots,t_{2n})}}(3)|
$
and
$
\mathcal U_{2n}=\langle\beta_0,\ldots,\beta_{2n}\rangle\subset |\mathcal O_{\mathbb P^{2n}_{(u_0,\ldots,u_{2n})}}(2)|.
$
Then
\begin{equation}
\label{TU-linear-systems}
\mathcal T_{2n}
=
|\mathcal I_{T_1,T_2,2T_3,T_4}(3)|,
\qquad
\mathcal U_{2n}
=
|\mathcal I_{U_1,U_2,U_3\rightarrow U_{4,n-1}}(2)|.
\end{equation}
Moreover, the rational maps
$$
\alpha_{2n}\colon \mathbb P^{2n}_{(t_0,\ldots,t_{2n})}\dashrightarrow \mathbb P^{2n}_{(u_0,\ldots,u_{2n})},
\qquad
t\longmapsto [\alpha_0(t):\ldots:\alpha_{2n}(t)]
$$
and
$$
\beta_{2n}\colon \mathbb P^{2n}_{(u_0,\ldots,u_{2n})}\dashrightarrow \mathbb P^{2n}_{(t_0,\ldots,t_{2n})},
\qquad
u\longmapsto [\beta_0(u):\ldots:\beta_{2n}(u)]
$$
are birational and inverse to each other on dense open subsets.

If $n=1$, the linear system $\mathcal T_2$ reduces to the conic system
$$
\langle
-2(t_0^2+t_0t_1+t_1^2),
t_2(2t_0+t_1),
t_1t_2
\rangle,
$$
which induces the standard plane Cremona transformation centered at the three points $[(-1+\xi)/2:1:0]$, $[(-1-\xi)/2:1:0]$ and $[0:0:1]$.
\end{Proposition}

\begin{proof}
The equalities of linear systems follow by direct substitution in the defining equations of the schemes $T_i$ and $U_i$. The forms $\alpha_0,\ldots,\alpha_{2n}$ vanish along $T_1$, $T_2$, $T_4$ and with multiplicity at least $2$ along $T_3$, while the forms $\beta_0,\ldots,\beta_{2n}$ satisfy the stated incidence and tangent conditions along $U_1,U_2,U_3$ and $U_4$. Conversely, these conditions impose the expected number of independent linear conditions on cubics and quadrics, so the displayed spans exhaust the two linear systems.

To prove birationality, one checks the compositions. Substituting the forms \eqref{alphaAEPs} into \eqref{betaAEPs} gives
$
\beta_{2n}\circ\alpha_{2n}=[t_0:\ldots:t_{2n}]
$
up to a common non-zero factor on the open set where both maps are defined. Conversely, substituting \eqref{betaAEPs} into \eqref{alphaAEPs} gives
$
\alpha_{2n}\circ\beta_{2n}=[u_0:\ldots:u_{2n}]
$
up to a common non-zero factor on the corresponding open set. Hence the maps are birational inverses. For $n=1$ the displayed conic system has three non-collinear base points over $k(\xi)$ and is the standard quadratic Cremona transformation.
\end{proof}

\begin{Remark}
The two parametrizations are related by the Cremona transformation of Proposition~\ref{quadro-cubic-cremona}. More precisely, the map $\overline\varphi_{2n}$ depends on the polynomials $\overline A,\overline B$ in \eqref{projectivepolynomialsAB}, while $\overline\phi_{2n}$ depends on $M_{2n},N_{2n}$ in \eqref{recurrencednm}. The maps $\alpha_{2n}$ and $\beta_{2n}$ identify these two descriptions on dense open subsets.
\end{Remark}

\begin{center}
\begin{tikzpicture}[baseline=(current bounding box.center), scale=1.2]
\node (Ptop) at (0,2) {$\mathbb P^{2n}_{(t_0,\ldots,t_{2n})}$};
\node (Pbot) at (0,0) {$\mathbb P^{2n}_{(u_0,\ldots,u_{2n})}$};
\node (X) at (4,1) {$X^{2n}_d$};

\draw[->, dashed] (Ptop) -- node[above right] {$\overline\phi_{2n},\ \deg 4d+4$} (X);
\draw[->, dashed] (Pbot) -- node[below right] {$\overline\varphi_{2n},\ \deg 2d+2$} (X);

\draw[->, dashed] (Ptop) to[bend left=25] node[right] {$\alpha_{2n},\ \deg 3$} (Pbot);
\draw[->, dashed] (Pbot) to[bend left=25] node[left] {$\beta_{2n},\ \deg 2$} (Ptop);
\end{tikzpicture}
\end{center}

\begin{Remark}
For $n=1$, the corresponding diagram becomes:
\end{Remark}

\begin{center}
\begin{tikzpicture}[baseline=(current bounding box.center), scale=1.2]
\node (Ptop) at (0,2) {$\mathbb P^{2}_{(t_0,t_1,t_2)}$};
\node (Pbot) at (0,0) {$\mathbb P^{2}_{(u_0,u_1,u_2)}$};
\node (X) at (4,1) {$X^2_d$};

\draw[->, dashed] (Ptop) -- node[above right] {$\overline\phi_2,\ \deg 2d+2$} (X);
\draw[->, dashed] (Pbot) -- node[below right] {$\overline\varphi_2,\ \deg 2d+2$} (X);

\draw[->, dashed] (Ptop) to[bend left=25] node[right] {$\alpha_2,\ \deg 2$} (Pbot);
\draw[->, dashed] (Pbot) to[bend left=25] node[left] {$\beta_2,\ \deg 2$} (Ptop);
\end{tikzpicture}
\end{center}
		
\section{Rational points over finite fields}
\label{Section_RatPoints}

We now analyze point counts over finite fields by using the explicit birational maps constructed above. We first treat the surface $X_d^2$.

For $2n=2$, the polynomials in \eqref{projectivepolynomialsAB} are
$$
\overline A
=
u_0^{2d+1}+(u_1+3u_2)(u_1^2+3u_2^2)^d,
\qquad
\overline B
=
u_0^{2d+1}+(u_1-u_2)(u_1^2+3u_2^2)^d.
$$
Let $Y^0=\{\overline A=\overline B=0\}\subset \mathbb P^2$. It has length $(2d+1)^2$ counted with multiplicities. When $d=1$, the scheme $Y^0$ has the following components:
$$
Y_1^0=\{u_0=u_1^2+3u_2^2=0\},
\qquad
Y_2^0=\{u_2=u_0^2-u_0u_1+u_1^2=0\},
\qquad
Y_3^0=\{u_2=u_0+u_1=0\}.
$$
Here $Y_1^0$ and $Y_2^0$ consist of two pairs of conjugate points over a quadratic extension of the base field, while $Y_3^0$ is one point over the base field. By Remark~\ref{osservazionetangentcone}, the two conjugate points in $Y_1^0$ have multiplicity $3$. These points are
$$
p_{1,+}=[0:\xi:1],
\quad
p_{2,+}=[1+\xi:2:0],
\quad
p_3=[1:-1:0],
\qquad
p_{1,-}=[0:-\xi:1],
\quad
p_{2,-}=[1-\xi:2:0],
$$
where $\xi^2=-3$.

For general $d$, the length of $Y^0$ is $4d^2+4d+1$. The components $Y_1^0$ and $Y_3^0$ do not depend on $d$, while $p_{1,+}$ and $p_{1,-}$ have multiplicity $2d^2+d$. The remaining component is
$$
Y_2^{0,d}
=
\left\{
u_2=
\sum_{i=0}^{2d}(-1)^iu_0^{2d-i}u_1^i
=
0
\right\},
$$
which consists of the $2d$ reduced points $p_{2,i}=[1:-\tau^i:0]$, where $\tau$ is a primitive $(2d+1)$-st root of unity and $i=1,\ldots,2d$.

We shall use the following elementary facts.

\begin{Lemma}
\label{x^2+3 lemma}
Let $q=p^m$. The polynomial $x^2+3\in \mathbb F_q[x]$ has a root in $\mathbb F_q$ if and only if $m$ is even, or $p\equiv 1\pmod 6$, or $p\in\{2,3\}$.
\end{Lemma}

\begin{Lemma}
\label{lemma of base field of 2d+1 points}
The polynomial $x^{2d+1}+1\in \mathbb F_q[x]$ has exactly $s=\gcd(q-1,2d+1)$ roots in $\mathbb F_q$.
\end{Lemma}

\begin{Corollary}
\label{numberpoints}
Let $X_d^2\subset \mathbb P^3$ be the surface over $\mathbb F_q$, where $q=p^m$, defined by
$$
(x_0+x_1)(x_0^2-x_0x_1+x_1^2)^d
+
(x_2+x_3)(x_2^2-x_2x_3+x_3^2)^d
=
0.
$$
If $p\neq 2,3$, set $s=\gcd(q-1,2d+1)$. Then
$$
|X_d^2(\mathbb F_q)|
=
\begin{cases}
q^2+sq+1, & q\equiv 5\pmod 6,\\
q^2+(s+4)q+1, & q\equiv 1\pmod 6.
\end{cases}
$$
In characteristic $3$, write $2d+1=3^ig$, with $3\nmid g$, and set $s'=\gcd(g,q-1)$. Then
$$
|X_d^2(\mathbb F_q)|
=
s'q^2+q+1.
$$
In characteristic $2$, the same formula holds with $q^2+sq+1$ when $m$ is odd and $q^2+(s+4)q+1$ when $m$ is even.
\end{Corollary}

\begin{Proposition}
\label{surface-blowup-model}
Assume $\operatorname{char}(k)\neq 2,3$. For $n=1$, the birational map $\overline\varphi_2\colon\mathbb P^2\dashrightarrow X_d^2$ is resolved by blowing up the base scheme of the plane linear system described below. More precisely, after blowing up the conjugate points $p_{1,+},p_{1,-}$ with the prescribed infinitely near tangent conditions and the simple base points $p_{2,i},p_3$, the induced morphism is birational onto $X_d^2$ and contracts precisely the curves corresponding to the exceptional and boundary components appearing in the Cremona model.
\end{Proposition}

\begin{proof}
The coordinates of $\overline\varphi_2$ are the generators of the plane linear system $\mathcal S_{2d+2,u}$ introduced below. The base conditions impose multiplicity $d+1$ at $p_{1,+}$ and $p_{1,-}$, with fixed principal tangents $L_{1,+}$ and $L_{1,-}$ of multiplicity $d$, and simple conditions at the remaining points. Blowing up this cluster removes the base scheme of the linear system. The resulting morphism agrees with $\overline\varphi_2$ outside the exceptional locus. Since Lemma~\ref{cremonaequivalence} identifies this model with the inverse quadratic system on $X_d^2$, the morphism is birational. This is the blow-up description used in Corollary~\ref{numberpoints}.
\end{proof}

Before proving Corollary~\ref{numberpoints}, we recall and adapt the plane Cremona description used in \cite{msslxa}. Assume $\operatorname{char}(k)\neq 2,3$. In $\mathbb P^2_{(u_0:u_1:u_2)}$ consider the points
$$
p_{1,+}=[0:\xi:1],
\qquad
p_{1,-}=[0:-\xi:1],
\qquad
p_{2,i}=[1:-\tau^i:0],
\qquad
p_3=[1:-1:0],
$$
where $i=1,\ldots,2d$ and $\tau$ is a primitive $(2d+1)$-st root of unity. Let $L_{1,+}=\{u_1-\xi u_2=0\}$ and $L_{1,-}=\{u_1+\xi u_2=0\}$. We denote by
$$
\mathcal S_{2d+2,u}
:=
\mathcal L^{2d+2}_{\substack{\xrightarrow{}\,dL_{1,\pm}\\ (d+1)p_{1,\pm},\ p_{2,i},\ p_3}}
$$
the linear system of plane curves of degree $2d+2$ with multiplicity at least $d+1$ at $p_{1,+}$ and $p_{1,-}$, with $L_{1,+}$ and $L_{1,-}$ as fixed principal tangents of multiplicity $d$, and passing through $p_{2,i}$ for $i=1,\ldots,2d$ and through $p_3$.

In $\mathbb P^2_{[v_0:v_1:v_2]}$, set
$$
q_{1,+}=[2\xi:-\xi:1],
\qquad
q_{1,-}=[-2\xi:\xi:1],
\qquad
q_{2,i}=[\tau^i-3:1:0],
\qquad
q=[1:0:0],
$$
and
$$
q_{3,+}=[2+4\xi:-\xi:1],
\qquad
q_{3,-}=[2-4\xi:\xi:1].
$$
Let
$
\mathcal S_{2d+1,v}:=\mathcal L^{2d+1}_{d(q_{1,\pm}),\,q_{2,i},\,d(q_{3,\pm})}
$
be the linear system of plane curves of degree $2d+1$ with multiplicity at least $d$ at $q_{1,+},q_{1,-},q_{3,+},q_{3,-}$ and passing through $q_{2,i}$ for $i=1,\ldots,2d$.

\begin{Lemma}
\label{cremonaequivalence}
The quadratic Cremona transformation
\begin{equation}
\label{cremona}
cr\colon \mathbb P^2_{[u_0:u_1:u_2]}\dashrightarrow \mathbb P^2_{[v_0:v_1:v_2]},
\qquad
[u_0:u_1:u_2]
\longmapsto
[3u_0^2+4u_0u_1+u_1^2+3u_2^2:-u_0^2-u_0u_1:u_0u_2]
\end{equation}
is induced by the conics through $p_{1,+},p_{1,-},p_3$ and gives a Cremona equivalence between $\mathcal S_{2d+2,u}$ and $\mathcal S_{2d+1,v}$.
\end{Lemma}

\begin{proof}
The inverse of $cr$ is
\begin{equation}
\label{cremonainverse}
cr^{-1}\colon \mathbb P^2_{[v_0:v_1:v_2]}\dashrightarrow \mathbb P^2_{[u_0:u_1:u_2]},
\qquad
[v_0:v_1:v_2]
\longmapsto
[v_1^2+3v_2^2:-v_0v_1-3v_1^2-3v_2^2:v_2(v_0+2v_1)].
\end{equation}
It is induced by the conics through $q_{1,+},q_{1,-}$ and $q$.

Let $C\in\mathcal S_{2d+2,u}$ be general and let $\Gamma=cr(C)$. Then
$$
\deg\Gamma
=
2\deg C-\operatorname{mult}_{p_{1,+}}C-\operatorname{mult}_{p_{1,-}}C-\operatorname{mult}_{p_3}C
=
4d+4-(d+1)-(d+1)-1
=
2d+1.
$$
The line $\langle p_{1,+},p_3\rangle$ is contracted to $q_{1,+}$. Since it meets $C$ in $2d+2$ points and $C$ has multiplicities $d+1$ at $p_{1,+}$ and $1$ at $p_3$, the image point $q_{1,+}$ has multiplicity $d$ on $\Gamma$. The same argument applies to $q_{1,-}$. The line $\langle p_{1,+},p_{1,-}\rangle$ is contracted to $q$, but it contributes no imposed multiplicity to $\Gamma$. The points $p_{2,i}$ are mapped bijectively to $q_{2,i}$, so $\operatorname{mult}_{q_{2,i}}\Gamma=1$.

Let $\widetilde{\mathbb P}^2_u$ be the blow-up of $\mathbb P^2_u$ at $p_{1,+},p_{1,-},p_3$, and let $\widetilde{\mathbb P}^2_v$ be the blow-up of $\mathbb P^2_v$ at $q_{1,+},q_{1,-},q$. The Cremona map lifts to an isomorphism $\widetilde{cr}\colon\widetilde{\mathbb P}^2_u\to\widetilde{\mathbb P}^2_v$:
$$
\begin{tikzcd}[row sep=large, column sep=large]
\widetilde{\mathbb P}^2_u
\arrow[r, bend left=20, "\widetilde{cr}"]
\arrow[d, "\pi_u"']
&
\widetilde{\mathbb P}^2_v
\arrow[l, bend left=20, "\widetilde{cr}^{-1}"']
\arrow[d, "\pi_v"]\\
\mathbb P^2_u
\arrow[r, bend left=20, dotted, "cr"]
&
\mathbb P^2_v
\arrow[l, bend left=20, dotted, "cr^{-1}"']
\end{tikzcd}
$$
Since $C$ has fixed principal tangents $L_{1,+}$ and $L_{1,-}$ at $p_{1,+}$ and $p_{1,-}$, both with multiplicity $d$, its strict transform meets the exceptional divisors over $p_{1,+}$ and $p_{1,-}$ in two fixed conjugate points. These are mapped to $q_{3,+}$ and $q_{3,-}$ by $\pi_v\circ\widetilde{cr}$. Thus $\Gamma\in\mathcal S_{2d+1,v}$.

The converse is identical. If $D\in\mathcal S_{2d+1,v}$ and $\Delta=cr^{-1}(D)$, then
$$
\deg\Delta
=
2\deg D-\operatorname{mult}_{q_{1,+}}D-\operatorname{mult}_{q_{1,-}}D-\operatorname{mult}_qD
=
4d+2-d-d
=
2d+2.
$$
The points $q_{1,+},q_{3,+},q$ are collinear, and their line is contracted to $p_{1,+}$; hence $\operatorname{mult}_{p_{1,+}}\Delta=d+1$, with fixed tangent $L_{1,+}$ of multiplicity $d$. The same computation applies to $p_{1,-}$. Finally, $q_{2,i}$ and $p_{2,i}$ correspond under the Cremona maps. Therefore $\Delta\in\mathcal S_{2d+2,u}$, and the two systems are Cremona equivalent.
\end{proof}

By Proposition~\ref{propolinsubsystem}, after a change of basis the system $\mathcal F_2$ is contained in $\mathcal S_{2d+2,u}$.

\begin{proof}[Proof of Corollary~\ref{numberpoints}]
By Lemma~\ref{cremonaequivalence}, the generators $\varphi_j\in\mathcal S_{2d+2,u}$, for $j=0,\ldots,3$, correspond to generators $\chi_j\in\mathcal S_{2d+1,v}=\mathcal L^{2d+1}_{d(q_{1,\pm}),\,q_{2,i},\,d(q_{3,\pm})}$. Let $E_{1,\pm}$, $E_{3,\pm}$ and $E_{2,i}$ be the exceptional divisors over $q_{1,\pm}$, $q_{3,\pm}$ and $q_{2,i}$, respectively.

Assume first that $p\neq 2,3$. If $q\equiv 1\pmod 6$, then $-3$ is a square in $\mathbb F_q$, so $\xi\in\mathbb F_q$ and the four points $q_{1,\pm},q_{3,\pm}$ are defined over $\mathbb F_q$. Hence the four exceptional divisors $E_{1,\pm}$ and $E_{3,\pm}$ are also defined over $\mathbb F_q$. By Lemma~\ref{lemma of base field of 2d+1 points}, exactly $s-1$ among the points $q_{2,i}$ are defined over $\mathbb F_q$; the missing root corresponds to the point not belonging to the base locus of the system. Therefore
$$
|X_d^2(\mathbb F_q)|
=
(|\mathbb P^2(\mathbb F_q)|-(s+3))+(s+3)|\mathbb P^1(\mathbb F_q)|
=
q^2+(s+4)q+1.
$$
If $q\equiv 5\pmod 6$, then the conjugate pairs $q_{1,\pm}$ and $q_{3,\pm}$ are not defined over $\mathbb F_q$, and the same count gives
$
|X_d^2(\mathbb F_q)|=q^2+sq+1.
$

If $p=3$, then
$$
X_d^2=(x_0+x_1)^{2d+1}+(x_2+x_3)^{2d+1}.
$$
Set $S=x_0+x_1$, $T=x_2+x_3$, and write $2d+1=g3^r$, with $3\nmid g$. Then
$
X_d^2=(S^g+T^g)^{3^r}.
$
Thus $X_d^2$ is the union, with multiplicities, of the planes $S=-\zeta T$ where $\zeta^g=1$. Exactly $s'=\gcd(g,q-1)$ of these planes are defined over $\mathbb F_q$, and they all contain the line $S=T=0$. Therefore
$$
|X_d^2(\mathbb F_q)|
=
s'(q^2+q+1)-(s'-1)(q+1)
=
s'q^2+q+1.
$$

Finally, in characteristic $2$ one uses the characteristic-$2$ Cremona model from Section~\ref{sec:char-two}. The same argument shows that the conjugate exceptional divisors are defined over $\mathbb F_{2^m}$ precisely when $m$ is even. Hence
$$
|X_d^2(\mathbb F_q)|
=
\begin{cases}
q^2+sq+1, & m\text{ odd},\\
q^2+(s+4)q+1, & m\text{ even}.
\end{cases}
$$
This completes the proof.
\end{proof}

We now record a simple point-counting lemma.

\begin{Lemma}
\label{hypersurfaceprojection}
Let $X\subset \mathbb P^{N+1}$ be the hypersurface over $\mathbb F_q$ defined by
$$
X=\{f(x_0,\ldots,x_N)+ax_{N+1}^D=0\},
$$
where $a\in\mathbb F_q^*$ and $\gcd(D,q-1)=1$. Then the projection $X\dashrightarrow \mathbb P^N$ induces a bijection $X(\mathbb F_q)\to\mathbb P^N(\mathbb F_q)$.
\end{Lemma}

\begin{proof}
The projection sends $[x_0:\ldots:x_N:x_{N+1}]$ to $[x_0:\ldots:x_N]$. Fix $y=[y_0:\ldots:y_N]\in\mathbb P^N(\mathbb F_q)$. Since $\gcd(D,q-1)=1$, the map $\mathbb F_q\to\mathbb F_q$, $c\mapsto c^D$, is bijective. Hence there is a unique $c\in\mathbb F_q$ such that $ac^D=-f(y_0,\ldots,y_N)$. Thus there is a unique point $[y_0:\ldots:y_N:c]\in X(\mathbb F_q)$ lying over $y$.
\end{proof}

\begin{Proposition}
\label{pointcountX2n}
Let $X_d^{2n}\subset\mathbb P^{2n+1}$ be defined by \eqref{X2n_equation} over $\mathbb F_q$. Assume $q\equiv 5\pmod 6$ and $\gcd(2d+1,q-1)=1$. Then
$$
|X_d^{2n}(\mathbb F_q)|
=
|\mathbb P^{2n}(\mathbb F_q)|.
$$
\end{Proposition}

\begin{proof}
The case $n=1$ follows from Corollary~\ref{numberpoints}. Assume the statement true for $n-1$. Write
$$
X_d^{2n}
=
\left\{
\sum_{i=0}^{n}
(x_{2i}+x_{2i+1})
(x_{2i}^2-x_{2i}x_{2i+1}+x_{2i+1}^2)^d
=
0
\right\}.
$$
For $[a_0:a_1]\in\mathbb P^1(\mathbb F_q)$, let $H_{[a_0:a_1]}=\{a_0x_{2n}+a_1x_{2n+1}=0\}$. If $[a_0:a_1]\neq [a'_0:a'_1]$, then
$
X_d^{2n}\cap H_{[a_0:a_1]}\cap H_{[a'_0:a'_1]}=X_d^{2n-2}.
$

Since $q\equiv 5\pmod 6$, the quadratic form $x_{2n}^2-x_{2n}x_{2n+1}+x_{2n+1}^2$ is anisotropic over $\mathbb F_q$. If $[a_0:a_1]\neq [1:-1]$, then
$$
H_{[a_0:a_1]}\cap X_d^{2n}
=
\{f(x_0,\ldots,x_{2n-1})+b x_{2n}^{2d+1}=0\}
$$
for some $b\in\mathbb F_q^*$, after choosing a coordinate on $H_{[a_0:a_1]}$. By Lemma~\ref{hypersurfaceprojection}, this hyperplane section has $|\mathbb P^{2n-1}(\mathbb F_q)|$ rational points. The points already lying in the common intersection $X_d^{2n-2}$ account, by the induction hypothesis, for $|\mathbb P^{2n-2}(\mathbb F_q)|$ points. Thus the contribution outside $X_d^{2n-2}$ is
$
|\mathbb P^{2n-1}(\mathbb F_q)|-|\mathbb P^{2n-2}(\mathbb F_q)|=q^{2n-1}.
$
There are $q$ such hyperplanes.

For $[a_0:a_1]=[1:-1]$, the section is the cone over $X_d^{2n-2}$ with vertex $P=[0:\cdots:0:1:-1]$. Hence it contributes
$$
q|X_d^{2n-2}(\mathbb F_q)|+1
=
q|\mathbb P^{2n-2}(\mathbb F_q)|+1
=
|\mathbb P^{2n-1}(\mathbb F_q)|.
$$
Therefore
$$
|X_d^{2n}(\mathbb F_q)|
=
q\cdot q^{2n-1}
+
|\mathbb P^{2n-1}(\mathbb F_q)|
=
q^{2n}
+
\frac{q^{2n}-1}{q-1}
=
\frac{q^{2n+1}-1}{q-1}.
$$
This is $|\mathbb P^{2n}(\mathbb F_q)|$.
\end{proof}

\begin{Remark}
For $H_0^{2n-1}=\{u_0=0\}\subset\mathbb P^{2n}$, the map $\overline\varphi_{2n}$ sends $H_0^{2n-1}$ to $X_d^{2n-2}$. Thus this hyperplane is contracted onto $X_d^{2n-2}$.
\end{Remark}

\begin{Proposition}
\label{pointcountY2n-2}
Let $k=\mathbb F_q$, with $q\equiv 5\pmod 6$, $\gcd(2d+1,q-1)=1$, and $n\geq 2$. Then
$$
|Y^{2n-2}(k)|
=
|\mathbb P^{2n-2}(k)|
=
\frac{q^{2n-1}-1}{q-1}.
$$
\end{Proposition}

\begin{proof}
Let $W^{2n-3}=Y^{2n-2}\cap H_0^{2n-1}$. Projection from the point $p=[0:\cdots:0:1]\in H_0^{2n-1}$, which does not belong to $Y^{2n-2}$, maps $W^{2n-3}$ to
$$
\overline W^{2n-3}
=
\left\{
\sum_{i=1}^{n-1}
v_{2i}(v_{2i-1}^2+3v_{2i}^2)^d
=
0
\right\}
\subset
\mathbb P^{2n-2}_{(v_0:\ldots:v_{2n-2})}.
$$
For a point $q=[q_0:\ldots:q_{2n-2}]\in\overline W^{2n-3}$, choose the first non-zero coordinate $q_i$. The fiber of $\pi_p|_{W^{2n-3}}$ over $q$ is defined on the line $\langle p,q\rangle$ by
$
q_i^{2d+1}u_{2n-2-i}^{2d+1}+u_{2n-1}^{2d+1}=0.
$
Since $\gcd(2d+1,q-1)=1$, this equation has exactly one solution over $k$. Thus $|W^{2n-3}(k)|=|\overline W^{2n-3}(k)|$. By the same argument as in \cite[Proposition~8.5]{msslxa}, one gets
$$
|W^{2n-3}(k)|
=
|\mathbb P^{2n-3}(k)|.
$$

The maps $\overline\varphi_{2n}$ and $\Theta_{2n}$ are inverse isomorphisms on the complements of their exceptional loci. Hence
$$
|\mathbb P^{2n}(k)|
-
|Y^{2n-2}(k)|
-
|H_0^{2n-1}(k)|
+
|W^{2n-3}(k)|
$$
equals
$$
|X_d^{2n}(k)|
-
|X_d^{2n-2}(k)|
-
\big(|Y^{2n-2}(k)|-|W^{2n-3}(k)|\big)|\mathbb P^1(k)|.
$$
Using Proposition~\ref{pointcountX2n} for $X_d^{2n}$ and $X_d^{2n-2}$, and using $|W^{2n-3}(k)|=|\mathbb P^{2n-3}(k)|$, this gives
$$
q|Y^{2n-2}(k)|
=
q\frac{q^{2n-1}-1}{q-1}.
$$
Therefore $|Y^{2n-2}(k)|=|\mathbb P^{2n-2}(k)|$.
\end{proof}

\begin{Remark}
Let $F_i=2^{2^i}+1$ be a Fermat prime. Since
$
\gcd(2d+1,F_i-1)=\gcd(2d+1,2^{2^i})=1,
$
Propositions~\ref{pointcountX2n} and~\ref{pointcountY2n-2} imply that, for $k=\mathbb F_{F_i}$,
$$
|X_d^{2n}(k)|=|\mathbb P^{2n}(k)|,
\qquad
|Y^{2n-2}(k)|=|\mathbb P^{2n-2}(k)|.
$$
Moreover, by Lemma~\ref{hypersurfaceprojection}, the hypersurfaces $\{\overline A=0\}$ and $\{\overline B=0\}$ have $|\mathbb P^{2n-1}(k)|$ rational points.
\end{Remark}

We also record a related generalization of Fermat hypersurfaces. Let $d$ be an odd prime and let $\delta\geq 1$. Define
$$
X_{d,\delta}^{2n}
=
\left\{
\sum_{i=0}^n
(x_{2i}+x_{2i+1})
\left(
\sum_{j=0}^{d-1}(-1)^j x_{2i}^{d-1-j}x_{2i+1}^j
\right)^\delta
=
0
\right\}
\subset
\mathbb P^{2n+1}.
$$
Its degree is $\Delta=\delta(d-1)+1$. For $\delta=1$, this is the Fermat hypersurface of degree $d$,
$$
X_{d,1}^{2n}
=
\{x_0^d+x_1^d+\cdots+x_{2n}^d+x_{2n+1}^d=0\}
\subset
\mathbb P^{2n+1}.
$$

\begin{Corollary}
\label{hypersurfaceprojection_cor}
If $k=\mathbb F_q$ and $\gcd(q-1,d)=1$, then
$$
|X_{d,1}^{2n}(k)|
=
|\mathbb P^{2n}(k)|.
$$
\end{Corollary}

\begin{Proposition}
\label{generalized-fermat-count}
Let $k=\mathbb F_q$ and $\delta\geq 1$. Assume that
$
\gcd(d,q-1)=\gcd(\delta(d-1)+1,q-1)=1.
$
Then
$$
|X_{d,\delta}^{2n}(k)|
=
|\mathbb P^{2n}(k)|.
$$
\end{Proposition}

\begin{proof}
The proof follows the same induction used in Proposition~\ref{pointcountX2n}. The condition $\gcd(d,q-1)=1$ gives the base case for $n=0$, since $X_{d,1}^0$ has the unique $k$-point $[1:-1]$. The condition $\gcd(\delta(d-1)+1,q-1)=1$ ensures that, at each induction step, the projection argument of Lemma~\ref{hypersurfaceprojection} applies to the last pair of variables.
\end{proof}

\section{Heights of rational points}
\label{sec:heights}

As an application of the rationality construction, we obtain lower bounds for the number of rational points of bounded height on $X_d^{2n}$.

Let $p\in\mathbb P^{2n+1}(\mathbb Q)$. Its reduced representative is the point $q=[q_0:\ldots:q_{2n+1}]$, with $q_i\in\mathbb Z$ and $\gcd(q_0,\ldots,q_{2n+1})=1$, such that $q=\lambda p$ for some $\lambda\in\mathbb Q^*$.

\begin{Definition}
The height of $p$ is
$
H(p)=\max\{|q_0|,\ldots,|q_{2n+1}|\}.
$
\end{Definition}

This definition extends to any number field; see \cite[Definition~1.5.4]{BG06}.

\begin{thm}
\label{height-lower-bound}
Let $K$ be a number field, and let $X^{2n}_{d,B}(K)$ be the set of $K$-rational points of $X_d^{2n}$ of height at most $B$ with respect to $\mathcal O_{X_d^{2n}}(1)$. Then, as $B\to\infty$,
$$
\#X^{2n}_{d,B}(K)
\gg
B^{\frac{2n+1}{2d+2}}.
$$
If $n=1$, then the lower bound improves to
$$
\#X^2_{d,B}(K)
\gg
B^{\frac{3}{2d+1}}.
$$
\end{thm}

\begin{proof}
The birational parametrization $\overline\varphi_{2n}\colon\mathbb P^{2n}\dashrightarrow X_d^{2n}$ in \eqref{projectivephis} is given by homogeneous polynomials of degree $2d+2$. Hence points of height at most $cB^{1/(2d+2)}$ in a non-empty open subset of $\mathbb P^{2n}$ are mapped to points of height at most $B$ on $X_d^{2n}$, for a constant $c>0$ depending only on the chosen height and on the coefficients of the parametrization.

Let $U\subset\mathbb P^{2n}$ be an open subset on which $\overline\varphi_{2n}$ is defined and generically finite onto its image. The complement has codimension at least one, and the standard height estimate on projective space, for instance \cite[Theorem~2.1]{Pey02}, gives
$$
\#\{u\in U(K):H(u)\leq B^{1/(2d+2)}\}
\gg
B^{\frac{2n+1}{2d+2}}.
$$
Since the map is generically finite on $U$, this gives the asserted lower bound.

For $n=1$, the Cremona equivalence in Lemma~\ref{cremonaequivalence} gives a parametrization of degree $2d+1$. Applying the same argument with degree $2d+1$ gives the exponent $3/(2d+1)$.
\end{proof}

\section{Cox rings and toric varieties}
\label{sec:cox-toric}

We conclude by describing the toric model naturally associated with the two conjugate $n$-planes. This also explains how the strict transform of $X_d^{2n}$ is represented by a single Cox equation of multidegree $(2d+1,-d,-d)$.

Let $k$ be the base field, and let $k(\xi)/k$ be the quadratic extension with $\xi^2=-3$. Recall that the two skew conjugate $n$-planes are defined over $k(\xi)$ by
$$
H_+
=
\{x_{2i}-a_+x_{2i+1}=0\text{ for }i=0,\ldots,n\},
\qquad
H_-
=
\{x_{2i}-a_-x_{2i+1}=0\text{ for }i=0,\ldots,n\},
$$
where $a_+=(1+\xi)/2$ and $a_-=(1-\xi)/2$.

We introduce homogeneous coordinates $y_0,\ldots,y_{2n+1}$ over $k(\xi)$ by the change of coordinates $\Phi$ given by
$$
y_{2j}=x_{2j}-a_+x_{2j+1},
\qquad
y_{2j+1}=x_{2j}-a_-x_{2j+1},
\qquad
j=0,\ldots,n.
$$
The inverse change of coordinates is
$$
x_{2j}
=
\frac{a_+y_{2j+1}-a_-y_{2j}}{\xi},
\qquad
x_{2j+1}
=
\frac{y_{2j+1}-y_{2j}}{\xi}.
$$
Under $\Phi$, the planes $H_+$ and $H_-$ become the coordinate subspaces
$$
H'_+
=
\{y_{2j}=0\text{ for }j=0,\ldots,n\},
\qquad
H'_-
=
\{y_{2j+1}=0\text{ for }j=0,\ldots,n\}.
$$
Thus, over $k(\xi)$, the blow-up $\mathcal T=\operatorname{Bl}_{H_+\cup H_-}\mathbb P^{2n+1}$ is isomorphic to
$
\mathcal T'=\operatorname{Bl}_{H'_+\cup H'_-}\mathbb P^{2n+1}.
$
Since $\mathcal T'$ is the blow-up of projective space along two coordinate linear subspaces, it is a smooth toric variety. Therefore $\mathcal T$ is a toric variety over $k(\xi)$, with descent datum induced by the Galois action.

The toric structure on $\mathcal T'$ is induced by the diagonal torus action on the $y$-coordinates,
$$
(t_0,\ldots,t_{2n+1})\cdot [y_0:\ldots:y_{2n+1}]
=
[t_0y_0:\ldots:t_{2n+1}y_{2n+1}].
$$
In the original $x$-coordinates, this action is obtained by conjugating by $\Phi$: if $g$ is the diagonal matrix representing the action in the $y$-coordinates, then its action in the $x$-coordinates is represented by $\Phi^{-1}g\Phi$. This torus action is not diagonal in the $x$-coordinates, but it is diagonal after base change to $k(\xi)$.

\begin{Remark}
Since $\Phi$ is block diagonal with $2\times 2$ blocks, any diagonal action with equal weights on each conjugate pair, say $t_{2i}=t_{2i+1}$ for every $i$, descends to an action defined over $k$. In particular, for
$
g=(h_0,h_0,h_1,h_1,\ldots,h_n,h_n)
$
with $h_i\in k^*$, the conjugated action $\Phi^{-1}g\Phi$ is defined over $k$ and is equal to the same diagonal action in the $x$-coordinates.

For the hypersurface $X_d^{2n}$, after normalizing $h_0=1$, this action sends the equation to
$$
\sum_{i=0}^n
h_i^{-(2d+1)}
(x_{2i}+x_{2i+1})
(x_{2i}^2-x_{2i}x_{2i+1}+x_{2i+1}^2)^d
=
0.
$$
Thus the diagonal subtorus preserving the decomposition into conjugate pairs moves $X_d^{2n}$ inside the same family.
\end{Remark}

\begin{Remark}
The union $H'_+\cup H'_-$ is defined in the $y$-coordinates by the $(n+1)^2$ quadrics
$
y_{2i}y_{2j+1}=0
$
for $0\leq i,j\leq n$. These correspond to the quadrics defining $H_+\cup H_-$ in \eqref{conjugateplanes}. In particular, the quadrics
$
x_{2i}^2-x_{2i}x_{2i+1}+x_{2i+1}^2
$
are mapped, up to a non-zero scalar, to $y_{2i}y_{2i+1}$, while the mixed quadrics are mapped to linear combinations of
$
y_{2i}y_{2j+1}
$
and
$
y_{2i+1}y_{2j}
$
with coefficients in $k(\xi)$.
\end{Remark}

We now rewrite the equation of $X_d^{2n}$ in the $y$-coordinates. Using the inverse change of coordinates, one gets
$
x_{2i}^2-x_{2i}x_{2i+1}+x_{2i+1}^2
$
equal, up to a non-zero scalar, to $y_{2i}y_{2i+1}$. Moreover,
$$
x_{2i}+x_{2i+1}
=
-\frac{1}{\xi}\big((a_-+1)y_{2i}-(a_++1)y_{2i+1}\big).
$$
Since $a_\pm+1=(3\pm\xi)/2$, the equation of $X_d^{2n}$ becomes, up to a non-zero scalar,
$$
\sum_{i=0}^n
y_{2i}^d y_{2i+1}^d
\big(3(y_{2i}-y_{2i+1})-\xi(y_{2i}+y_{2i+1})\big)
=
0.
$$
Define
$$
S
=
\sum_{i=0}^n
y_{2i}^d y_{2i+1}^d(y_{2i}+y_{2i+1}),
\qquad
D
=
\sum_{i=0}^n
y_{2i}^d y_{2i+1}^d(y_{2i}-y_{2i+1}).
$$
Then the equation of $X_d^{2n}$ over $k(\xi)$ is
$
3D-\xi S=0.
$
Let $\sigma$ be the non-trivial automorphism of $k(\xi)/k$, so $\sigma(\xi)=-\xi$. Since $\sigma(y_{2i})=y_{2i+1}$ and $\sigma(y_{2i+1})=y_{2i}$, we have $\sigma(S)=S$ and $\sigma(D)=-D$. Hence $S$ and $D/\xi$ are defined over $k$ when written in the original $x$-coordinates. The single equation $3D-\xi S=0$ is therefore compatible with the Galois descent from $k(\xi)$ to $k$.

\begin{Remark}
More generally, consider a hypersurface of the form
$$
\sum_{i=0}^n
(a_{2i}x_{2i}+a_{2i+1}x_{2i+1})
(x_{2i}^2-x_{2i}x_{2i+1}+x_{2i+1}^2)^d
=
0,
$$
with $a_j\in k$. In the $y$-coordinates its equation can again be written as $3D_a-\xi S_a=0$, where
$$
S_a
=
\sum_{i=0}^n
y_{2i}^d y_{2i+1}^d
\left(\frac{a_{2i}}{2}\right)
(y_{2i}+y_{2i+1}),
$$
and
$$
D_a
=
\sum_{i=0}^n
y_{2i}^d y_{2i+1}^d
\left(-\frac{a_{2i}}{2}-a_{2i+1}\right)
(y_{2i}-y_{2i+1}).
$$
Again $S_a$ is Galois-invariant and $D_a$ is Galois anti-invariant, so $S_a$ and $D_a/\xi$ descend to polynomials over $k$ in the original coordinates.
\end{Remark}

The Cox ring of $\mathcal T'=\operatorname{Bl}_{H'_+\cup H'_-}\mathbb P^{2n+1}$ over $k(\xi)$ is
$
R=k(\xi)[z_0,\ldots,z_{2n+1},w_+,w_-],
$
where $z_i$ corresponds to the strict transform of the coordinate hyperplane $\{y_i=0\}$, and $w_+$, $w_-$ correspond to the exceptional divisors $E_+$ and $E_-$ over $H'_+$ and $H'_-$. Since we blow up two disjoint $n$-planes, $\rho(\mathcal T')=3$ and
$
\operatorname{Pic}(\mathcal T')=\mathbb ZH\oplus\mathbb ZE_+\oplus\mathbb ZE_-.
$
The $\mathbb Z^3$-grading is given by
$$
\deg(z_{2i})=(1,-1,0),
\qquad
\deg(z_{2i+1})=(1,0,-1),
\qquad
\deg(w_+)=(0,1,0),
\qquad
\deg(w_-)=(0,0,1).
$$
Equivalently, the grading is encoded by
$$
\begin{array}{c|cccc}
& z_0,z_2,\ldots,z_{2n} & z_1,z_3,\ldots,z_{2n+1} & w_+ & w_-\\
\hline
H & 1 & 1 & 0 & 0\\
E_+ & -1 & 0 & 1 & 0\\
E_- & 0 & -1 & 0 & 1
\end{array}
$$
The irrelevant ideal is
$$
\mathfrak B
=
(z_0,\ldots,z_{2n+1})
\cap
(w_+,z_1,z_3,\ldots,z_{2n+1})
\cap
(w_-,z_0,z_2,\ldots,z_{2n}),
$$
and
$
\mathcal T'=(\operatorname{Spec}R\setminus V(\mathfrak B))/(\mathbb G_m)^3.
$

Let $X\subset\mathbb P^{2n+1}$ be any hypersurface of degree $2d+1$ with multiplicity $d$ along both $H'_+$ and $H'_-$. Its total transform has class $(2d+1)H$, while its strict transform has class
$
(2d+1)H-dE_+-dE_-.
$
In the above $\mathbb Z^3$-grading this is the degree $(2d+1,-d,-d)$.

For our hypersurface, the strict transform $\widetilde X_d^{2n}\subset\mathcal T'_{k(\xi)}$ is cut out by the strict transform of the equation $3D-\xi S=0$. In Cox coordinates this equation is
\begin{equation}
\label{cox-equation-strict-transform}
3\widehat D-\xi\widehat S=0,
\end{equation}
where
$$
\widehat S
=
\sum_{i=0}^n
z_{2i}^d z_{2i+1}^d(w_+z_{2i}+w_-z_{2i+1}),
\qquad
\widehat D
=
\sum_{i=0}^n
z_{2i}^d z_{2i+1}^d(w_+z_{2i}-w_-z_{2i+1}).
$$
Both $\widehat S$ and $\widehat D$ have multidegree $(2d+1,-d,-d)$. Indeed, each monomial
$
z_{2i}^d z_{2i+1}^d w_+z_{2i}
$
has degree
$
d(1,-1,0)+d(1,0,-1)+(0,1,0)+(1,-1,0)=(2d+1,-d,-d),
$
and similarly for
$
z_{2i}^d z_{2i+1}^d w_-z_{2i+1}.
$

Thus $\widetilde X_d^{2n}$ is a divisor of class $(2d+1)H-dE_+-dE_-$ on $\mathcal T'_{k(\xi)}$. The Galois action exchanges $z_{2i}$ with $z_{2i+1}$ and $w_+$ with $w_-$; under this action $\widehat S$ is invariant and $\widehat D$ is anti-invariant. Consequently, the divisor defined by \eqref{cox-equation-strict-transform} is stable under the descent datum and gives the desired $k$-model of the strict transform on $\mathcal T$.

\bibliographystyle{alpha}
\bibliography{Biblio}

@article{msslxa,
  author  = {Massarenti, Alex},
  title   = {Rational points on even dimensional {F}ermat cubics},
  journal = {Transactions of the London Mathematical Society},
  volume  = {13},
  number  = {3},
  year    = {2026},
  pages   = {e70028},
  doi     = {10.1112/tlm3.70028},
  eprint  = {2406.07223},
  archivePrefix = {arXiv},
  primaryClass = {math.AG}
}

@book{BG06,
  author    = {Bombieri, Enrico and Gubler, Walter},
  title     = {Heights in Diophantine Geometry},
  series    = {New Mathematical Monographs},
  volume    = {4},
  publisher = {Cambridge University Press},
  year      = {2006}
}

@article{Pey02,
  author  = {Peyre, Emmanuel},
  title   = {Points de hauteur born{\'e}e et g{\'e}om{\'e}trie des vari{\'e}t{\'e}s (d'apr{\`e}s Y. Manin et al.)},
  journal = {Ast{\'e}risque},
  volume  = {282},
  year    = {2002},
  pages   = {323--344},
  note    = {S{\'e}minaire Bourbaki, exp. no. 891}
}

@article{BT98,
  author  = {Batyrev, Victor V. and Tschinkel, Yuri},
  title   = {Manin's conjecture for toric varieties},
  journal = {J. Algebraic Geom.},
  volume  = {7},
  number  = {1},
  year    = {1998},
  pages   = {15--53},
  eprint  = {alg-geom/9510014}
}

@article{MS20,
  author  = {McKinnon, David and Satriano, Matthew},
  title   = {Approximating rational points on toric varieties},
  journal = {Trans. Amer. Math. Soc.},
  volume  = {373},
  number  = {12},
  year    = {2020},
  pages   = {8471--8509},
  doi     = {10.1090/tran/8138},
  eprint  = {2004.05212}
}

@article{Pie24,
  author  = {Pieropan, Marta and Schindler, Damaris},
  title   = {Hyperbola method on toric varieties},
  journal = {J. {\'E}c. polytech. Math.},
  volume  = {11},
  year    = {2024},
  pages   = {107--157},
  doi     = {10.5802/jep.251},
  eprint  = {2001.09815}
}

@article{PS24,
  author  = {Pieropan, Marta and Schindler, Damaris},
  title   = {Points of bounded height on certain subvarieties of toric varieties},
  journal = {Algebra Number Theory},
  year    = {2025},
  volume  = {19},
  number  = {11},
  pages   = {2281--2306},
  doi     = {10.2140/ant.2025.19.2281},
  eprint  = {2403.19397}
}

@misc{Bon24,
  author       = {Bongiorno, Nicolas},
  title        = {Multi-height analysis of rational points of toric varieties},
  year         = {2024},
  eprint       = {2412.04226},
  archivePrefix = {arXiv},
  primaryClass = {math.NT}
}
	
\end{document}